\documentclass[english,9pt]{amsart}

\usepackage{multirow}
\usepackage{bm}
\usepackage{amsmath}
\usepackage{amssymb}
\usepackage{amsbsy}
\usepackage{amsfonts}
\usepackage{array}
\usepackage{verbatim}
\usepackage{graphicx}
\usepackage{listings}
\usepackage{subfigure}
\usepackage{multirow}
\usepackage[normalem]{ulem}
\usepackage{cancel}

\graphicspath{{./figures/},{./}}

\usepackage{url}

\usepackage{tikz}
\usetikzlibrary{matrix,arrows,decorations.pathmorphing}

\newcommand{\vect}[1]{\bm{#1}}

\newcommand{\jmp}[1]{\,[\![#1]\!]}

\definecolor{meogreen}{rgb}{0.3,0.6,0.1}

\newcommand{\unit}[1]{\mathrm{#1}}
\newcommand{\electronvolt}{\unit{eV}}
\newcommand{\kelvin}{\unit{K}}
\newcommand{\nano}{\unit{n}}
\newcommand{\meter}{\unit{m}}
\newcommand{\second}{\unit{s}}
\newcommand{\volt}{\unit{V}}

\begin{document}

\title{Multiscale Modeling and Simulation  of Organic Solar Cells}

\author{C. de Falco$^{1,2}$ \and M. Porro$^{1,3}$ \and R. Sacco$^{1}$ \and M. Verri$^{1}$}

\address{$^{1}$ Dipartimento di Matematica \lq\lq F. Brioschi\rq\rq,
               Politecnico di Milano, \\Piazza L. da Vinci 32, 20133 Milano, Italy\\
         $^{2}$      MOX Modeling and Scientific Computing\\
         $^{3}$      Center for Nano Science and Technology @PoliMi, \\
                   Istituto Italiano di Tecnologia,\\
               via Pascoli 70/3, 20133 Milano, Italy}

\date{\today}

\maketitle

\begin{abstract}
In this article, we continue our mathematical study
of organic solar cells (OSCs)
and propose a two-scale (micro- and macro-scale) 
model of heterojunction OSCs with interface geometries 
characterized by an arbitrarily complex morphology. 
The microscale model consists of a system of 
partial and ordinary differential
equations in an heterogeneous domain, that provides a full 
description of excitation/transport phenomena occurring 
in the bulk regions and dissociation/recombination
processes occurring in a thin material slab across the interface.
The macroscale model is obtained by a micro-to-macro scale 
transition that consists of averaging the mass balance 
equations in the normal direction across the interface thickness,
giving rise to nonlinear transmission conditions 
that are parametrized by the interfacial width.
These conditions account in a lumped manner 
for the volumetric dissociation/recombination phenomena
occurring in the thin slab and depend locally
on the electric field magnitude and orientation.
Using the macroscale model in two spatial dimensions, device 
structures with complex interface
morphologies, for which existing data are available, 
are numerically investigated showing that, 
if the electric field orientation
relative to the interface is taken into due account, the device
performance is determined not only by the total
interface length but also by its shape.

\end{abstract}


\section{Introduction and Motivation}
\label{sec:intro}
Research on photovoltaic energy conversion has recently received 
great impulse due to the growing demand for low carbon dioxide 
emission energy sources.
In particular, the high manufacturing cost of crystalline silicon 
and the latest advancements on semiconducting polymer design and 
synthesis in recent years have directed the attention of the scientific community towards Organic Solar Cells (OSCs), {\it i.e.}\ 
solar cells based on organic
materials~\cite{coakley_mcgehee,Gunes2007,halls,mcgehee,mihailetchi2004,pope1999,tang}, 
especially because of the very limited thermal 
budget required  for the production of such materials and of their 
amenability to be deposited on large areas, which is fundamental 
in light harvesting applications. 
One of the main peculiarities of OSCs is that
most physical phenomena that are critical for charge photogeneration
occur at the interface between the two materials that constitute the 
active layer of such devices.
In order to increase cell efficiencies, currently of the order of 
about 10\%~\cite{heliatek_record}, the optimization of the morphology
of such interface is considered by device designers to be an 
issue at least as important as the optimization of the donor and acceptor 
optoelectronic characteristics~\cite{Ray2012}.

For this reason, in this article we continue our mathematical study
of organic photovoltaic device models started in~\cite{defalco}, 
and we focus on the accurate and computationally efficient
modeling of the main dissociation/recombination
processes occurring in a thin material slab across the 
material interface and evaluating their impact on device photoconversion
performance. 
% Looking at existing literature, we can
% see that electrical transport models fall into three main categories:
% equivalent circuits, continuum models that solve the continuity
% and Poisson equations, and microscopic models based on
% Monte Carlo simulation. Most continuum models consider
% the whole OSC as a homogeneous medium and provide quantitatively 
% accurate comparisons with experiments~\cite{mihailetchi2007}.
% However, whilst computationally efficient, 
% continuum models are unable to describe
% carrier interactions in sufficient detail, 
% and they \nuovo{synthetize} the influence of morphology 
% \nuovo{into bulk parameters such as ?? QUALI SONO I PARAMETRI??
% SPECIFICARLI, DAREBBE IMPATTO ALLA PRESENTAZIONE, ORA E' TROPPO VAGO}. 
% To include \nuovo{carrier interaction effects}, 
% one should resort to a microscopical model based on dynamical
% Monte Carlo simulation, as done in~\cite{marsh2007}.
% However, as Coulomb interactions are hard to treat accurately in Monte
% Carlo simulations because of their long range nature~\cite{casalegno}, 
With this aim, we consider a two-scale approach to OSC 
simulation that is intermediate between 
a continuum model~\cite{mihailetchi2007} 
and a full microscopic model~\cite{marsh2007,casalegno},
% The macroscale model introduced in this article is 
% a mathematical rationale of that proposed 
%in~\cite{williams_th,williams2008}, 
and represents an extension 
to the case of arbitrary interface geometries of the one-dimensional model 
for bilayer OSC devices proposed in~\cite{barker2003}.

The approach is based on the introduction of two distinct 
levels of description of the physical system at hand, 
a micro and a macro scale, and of two corresponding mathematical
models based on classical mass balance conservation laws.
At the microscale, a system of PDEs in a heterogeneous domain 
provides a full description of the excitation/transport phenomena occurring 
in the bulk regions and of the dissociation/recombination 
processes occurring in a thin material slab across the interface.
%Despite being ready in principle for use in simulation purposes,
The numerical treatment of the microscale model presents several
difficulties related to the wide difference in size between the 
bulk regions and the interfacial width $H$. 
As a matter of fact, as polaron dissociation is assumed to 
occur in the first layer of polymer chains on either side 
of the interface, the length scale $H$ 
can be taken to be that of the average separation between polymer chains
which is typically more than two orders of magnitude smaller than
the size of the bulk regions~\cite{barker2003,buxton2007,williams2008}. 
% This value is
% quite small compared to the average size of the bulk regions, which
% is typically more than two orders oof the order of $100\,\nano\meter$.

Therefore, to obtain a computationally efficient model, we 
carry out a micro-to-macro scale transition that somewhat resembles
model-reduction techniques used for porous media 
with thin fractures~\cite{Jaffre2005}, for reaction problems 
with moving reaction fronts~\cite{transitions} and
for electrochemical transport across biological 
membranes~\cite{Mori2009}, and relies on a systematic 
averaging of the mass balance equations 
in the normal direction across the interface thickness.
The resulting macroscale model is a system of incompletely parabolic 
PDEs to describe mass transport in the materials, nonlinearly
coupled with ODEs and transmission  conditions localized 
at the heterojunction parametrized by the interfacial
width $H$. These conditions account in a lumped manner 
for the volumetric dissociation/recombination phenomena
occurring in the interfacial thin slab.
The fact that in the macroscale model the interface is reduced 
to a zero width surface is further exploited to account for the
local dependence of the polaron dissociation rate on 
the electric field orientation, which is the main advantage 
--together with the computational cost reduction-- of our
approach, as compared to previous multi-dimensional 
models~\cite{buxton2007,williams2008,jerome2010}.

%These conditions are  and depend locally on the electric field magnitude 
%and orientation,  
%so that, using the novel macroscale formulation, it is
%straightforward to include the dependence
%of the dissociation rate on the electric field and on 
%the morphology of the material interface.

An outline of the article is as follows.
In Sect.~\ref{sec:NBHJSCs} we illustrate the sequence
of physical phenomena that lead from photon absorbtion 
to current harvesting in an OSC.
Sect.~\ref{sec:math_model} is devoted to characterizing
the mathematical model of an OSC.
In Sect.~\ref{sec:geometry} 
we describe the geometrical heterogeneous structure of the 
device, while in Sect.~\ref{sec:basic_assumptions} 
we introduce the basic modeling assumptions 
on the dependent variables of the problem. Then, in 
Sect.~\ref{sec:microscale_model} we present
the microscale PDE/ODE model system with the initial 
and boundary conditions, while in Sect.~\ref{sec:macroscale_model} 
we describe in detail the scale transition procedure that
leads from the microscale model to the macroscale equation system.
We complete the mathematical picture of a bilayer OSC by illustrating in Sect.~\ref{sec:k_diss_model}
a novel model that we have devised for including the dependence 
of the polaron dissociation rate constant on 
the local electric field and on the morphology of the material
interface.
%discussing also the transmission conditions that arise from 
%the average procedure across the interface thickness.
In Sect.~\ref{sec:numer_approx} we briefly comment 
about the numerical methods used for discretizing
the macroscale model, while Sect.~\ref{sec:num_res}
is devoted to presenting and discussing numerical results.
In particular, in Sect.~\ref{sec:comparison_macro_micro} we successfully 
perform the validation of the accuracy of the macroscale model
with respect to the microscale model through the numerical
simulation of a one--dimensional OSC under different working
conditions. Extensive simulations of two--dimensional OSC structures
are instead reported in Sects.~\ref{sec:comparison} to
\ref{sec:morphology_2} in order to both validate the proposed 
macroscale model with
respect to previously available numerical results and to
analyze its effectiveness in the study of complex interface morphologies.
Finally, in Sect.~\ref{sec:concl} we draw some conclusions and sketch possible 
directions for further research in the area of OSC modeling and simulation.

\section{Basic Principles of Photocurrent Generation in OSCs}
\label{sec:NBHJSCs}
\begin{figure}[h!]
\centering
\subfigure[Basic configuration]{
\includegraphics[height=.15\textheight]{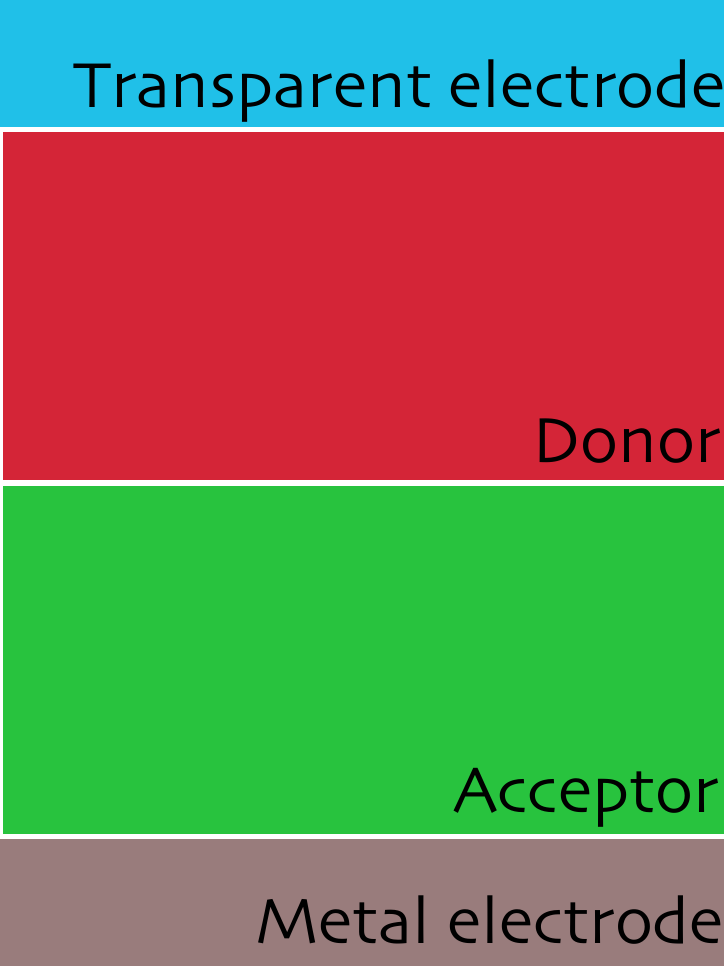}
\label{fig_osc_design}
}
\;
\subfigure[Device with disordered interface morphology]{
\includegraphics[height=.15\textheight]{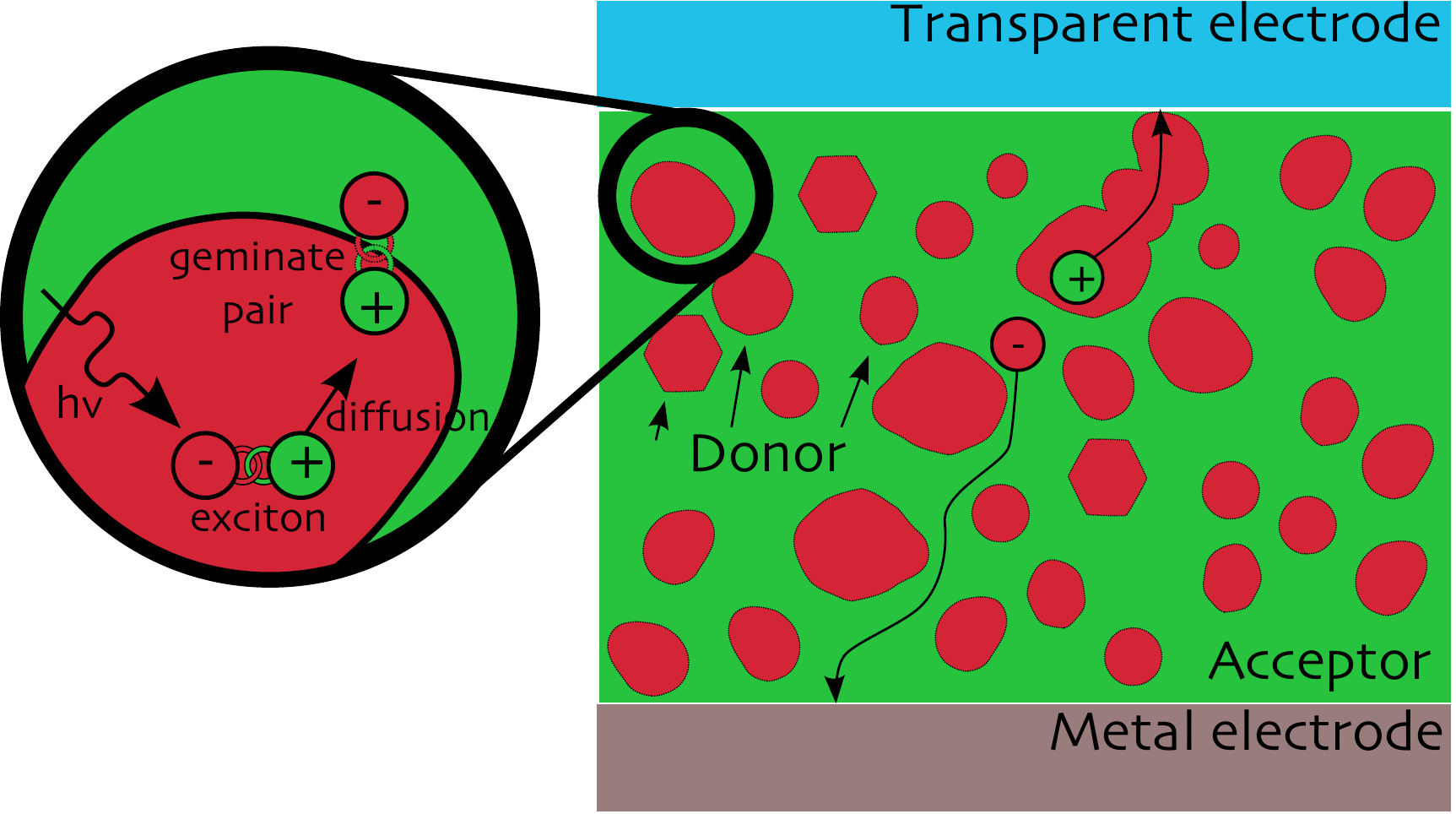}
\label{fig_bulk_heterojunction}
}
\caption{Structure of an organic solar cell.}
\label{fig:structure}
\end{figure}

In this section, we describe the basic principles of
photocurrent generation in OSCs only to the extent 
strictly needed for understanding the naming conventions adopted in 
the following sections. For a more thorough introduction to
the subject, we refer the interested reader 
to~\cite{Gunes2007,pope1999,mcgehee}.
The typical structure of an OSC is constituted by a thin film a cross-section
of which is schematically represented in Fig.~\ref{fig:structure}.

The photoactive layer of the device consists of two materials, 
one with higher electron 
affinity (the ``acceptor'', for example F8BT, P3HT) and one with lower electron affinity 
(the ``donor'', for example PFB, PCBM), 
sandwiched between two electrodes, one of which is transparent to allow light to
enter the photoactive layer while the other is reflecting in order to increase the light path through the device.

\begin{figure}[h!]
\centering
\includegraphics[width=0.9\textwidth]{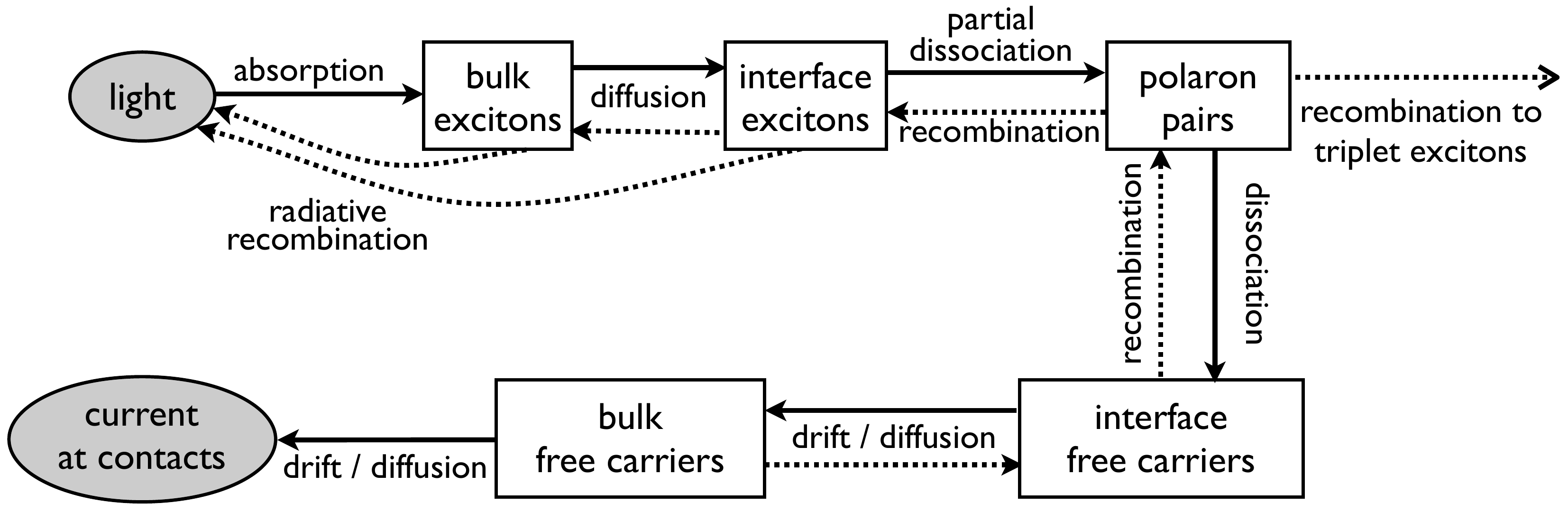}
\caption{Flow-chart of the photoconversion mechanisms in an OSC.}
\label{fig_flowchart}
\end{figure}

The sequence of physical phenomena that leads from photon absorption to current harvesting at
the device contacts is represented in Fig.~\ref{fig_flowchart}.

Absorption of a photon in either material produces an electron-hole pair, usually referred 
to as an \emph{exciton} whose binding energy is of the order of about $0.5 \div 1\,\electronvolt$.
Excitons may diffuse through the device until they either recombine or reach the interface 
between the donor and acceptor phases. 
If this latter event occurs, the exciton may get trapped at the interface in such a way that its electron 
component lays in the high electron affinity region while the hole component 
lays in the low electron affinity region.
Such a trapped excited state is referred to as a \emph{polaron pair} or \emph{geminate
pair}~\cite{barker2003,williams2008,morteani2008,mingebach2012}
and has a lower binding energy compared to that of the exciton state, as the Coulomb attraction between the electron and hole
is reduced by the chemical potential drop between the two materials.
The polaron binding force may be overcome by the electric field induced by the small built-in voltage between the metallic contacts
thus leading to the formation of two independent charged particles (one electron and one hole), 
otherwise the polaron pair may return to the untrapped exciton state
or recombine.
Free charge carriers move by drift and diffusion mechanisms and, unless 
they are captured along their path 
 by the coulombic attraction of an oppositely charged particle and recombine at the interface to form a new polaron pair, they eventually
reach the contacts thus producing a measurable external current.

%It is important to notice that the physical phenomena that are most crucial in determining the efficiency of the cell all occur at the interface between the acceptor and donor material.
%For this reason, the focus of the model developed in the present paper is to allow to accurately describe the interface morphology and to estimate its impact on device performance.

%%% Local Variables: 
%%% mode: latex
%%% TeX-master: "paper_osc_multiscale"
%%% End: 

\section{Mathematical Model}
\label{sec:math_model}
In this section, we propose a PDE/ODE model of photoconversion 
and charge transport mechanisms in an OSC. 
The model relies on a two-scale approach that is based on the introduction 
of two distinct levels of description of the physical system, 
namely, a micro and a macro scale, and of two corresponding 
mathematical equation systems based on classical mass 
balance conservation laws. The construction of the model 
proceeds through four steps. In Sect.~\ref{sec:geometry}, 
we describe the geometrical and heterogeneous structure of the device
which consists of two bulk regions (the acceptor and donor phases)
separated by an interface region of (finite) thickness $2H$,
while in Sect.~\ref{sec:basic_assumptions}, 
we introduce the basic modeling assumptions 
on the dependent variables of the problem.
In Sect.~\ref{sec:microscale_model}, we introduce 
the microscale PDE/ODE model system of conservation laws
that governs transport of the various species throughout the device,
together with its initial and boundary conditions
and the generation/recombination mechanisms that occur
in each subdomain of the heterogeneous device.
In Sect.~\ref{sec:macroscale_model}, we describe in detail 
the scale transition procedure that is applied to the microscale 
model in order to obtain the macroscale equation system.
This latter system basically consists of the same equations
as in the microscale model, but satisfied in the separate 
acceptor and donor phases, coupled through a set of 
flux transmission conditions across
the material interface that synthetize in a \lq\lq lumped\rq\rq{}
manner the dissociation and recombination mechanisms
that actually occur in the thin volumetric slab of width $2H$
surrounding the interface itself.
The resulting macroscale model system is a compromise between 
a continuum model and a full microscopic model, and
represents a consistent mathematical rationale and 
generalization of the various models
proposed in~\cite{barker2003,williams_th,williams2008}.
We conclude our mathematical picture of the OSC by illustrating 
in Sect.~\ref{sec:k_diss_model} a novel model 
of the polaron dissociation rate properly devised 
for including the dependence on the local electric field 
and on the morphology of the material interface.

\subsection{Geometry of the Heterogeneous Computational 
Domain}\label{sec:geometry}

\begin{figure}[h!]
\centering
\subfigure[3D OSC]
{\includegraphics[width=0.45\textwidth]{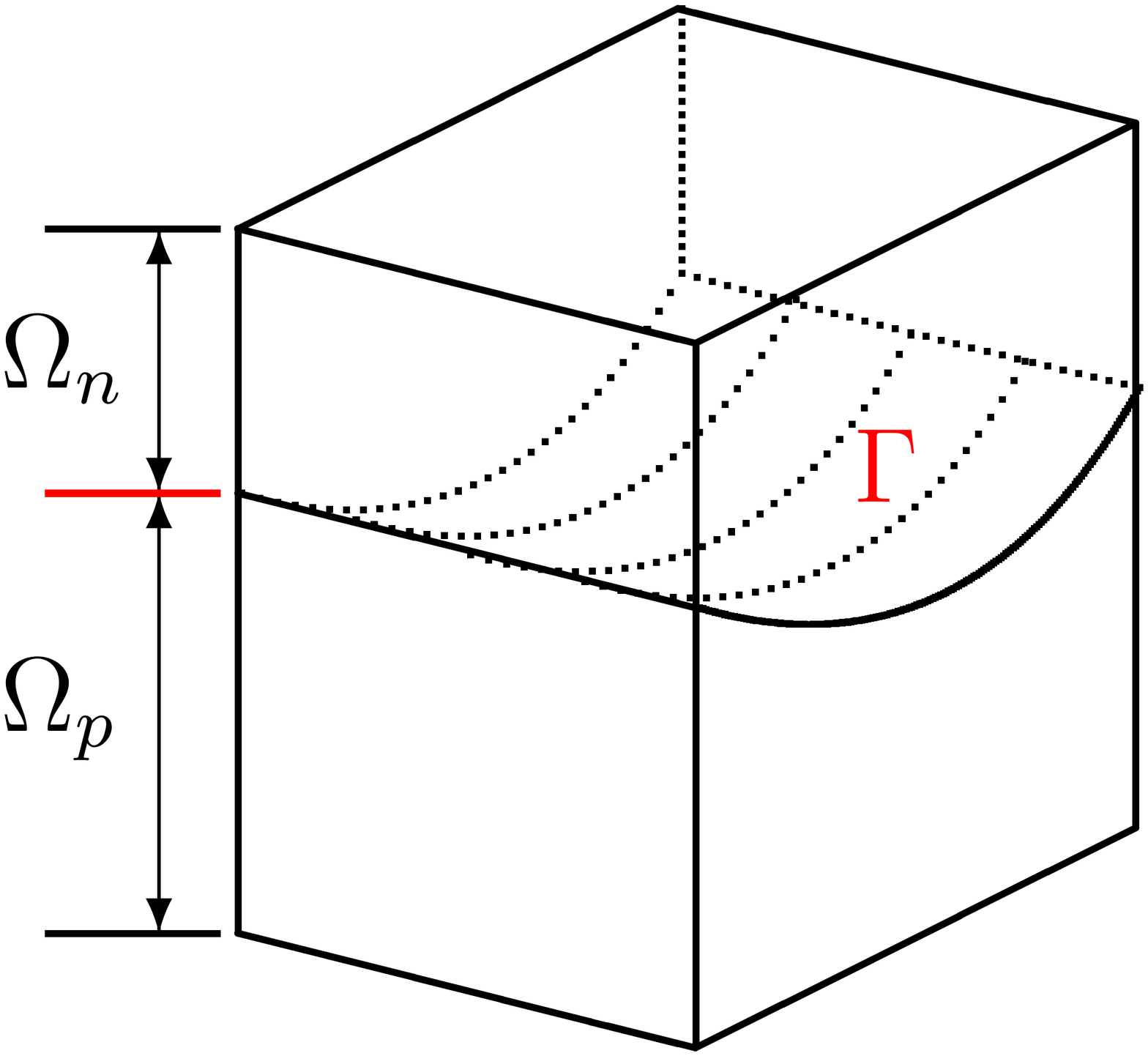}
\label{fig:Omega_3d}}
\subfigure[2D cross-section]
{\includegraphics[width=0.45\textwidth]{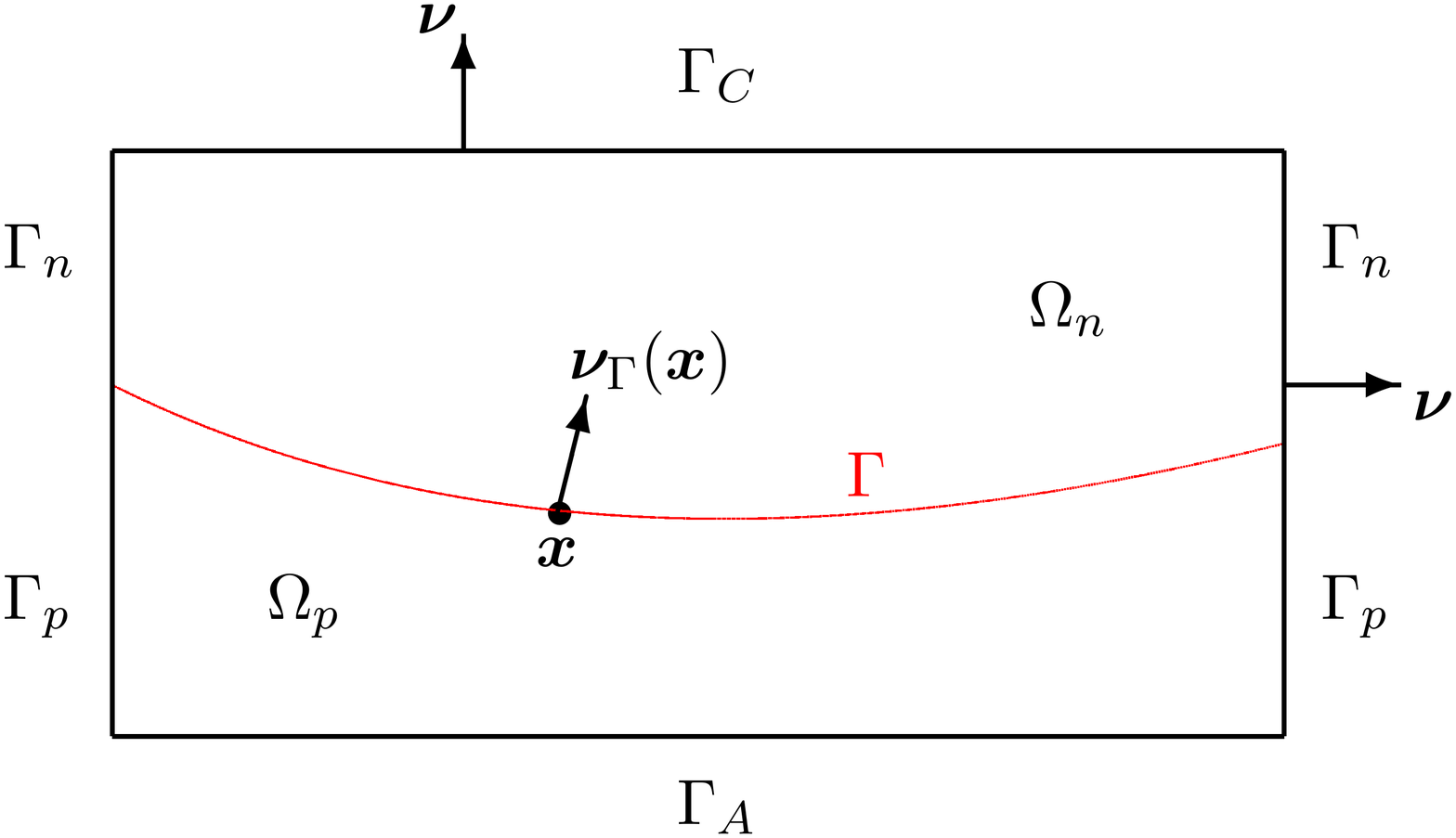}
\label{fig:Omega_2d}}
\caption{Geometry of the cell bulk.}
\label{fig:fig_Omega}
\end{figure}

A schematic 3D picture of the OSC is illustrated in
Fig.~\ref{fig:Omega_3d}. 
The device structure $\Omega$ is a parallelepiped-shaped open subset of 
$\mathbb{R}^3$ divided into two open disjoint subregions, 
$\Omega_n$ (acceptor) and $\Omega_p$ (donor), 
separated by a regular oriented surface $\Gamma = 
\partial \Omega_n \cap \partial \Omega_p$~\cite{dautraylionsvol2}
on which, for each $\vect{x} \in \Gamma$, we can define a unit normal
vector $\boldsymbol{\nu}_{\Gamma}(\vect{x})$
directed from $\Omega_p$ into $\Omega_n$.
The top and bottom surfaces of the structure 
are mathematical representations of
the cell electrodes, cathode and anode, denoted as $\Gamma_C$ and
$\Gamma_A$, respectively, in such a way that 
$\partial \Omega_n=\Gamma_C \cup \Gamma \cup \Gamma_n$ and
$\partial \Omega_p=\Gamma_A \cup \Gamma \cup \Gamma_p$ 
(see Fig.~\ref{fig:Omega_2d}). 
We also denote by $\boldsymbol{\nu}$ the 
unit outward normal vector over the cell boundary 
$\partial\Omega$.

\begin{figure}[h!]
\centering
\subfigure[3D representation of the interface]
{\includegraphics[width=0.45\textwidth]{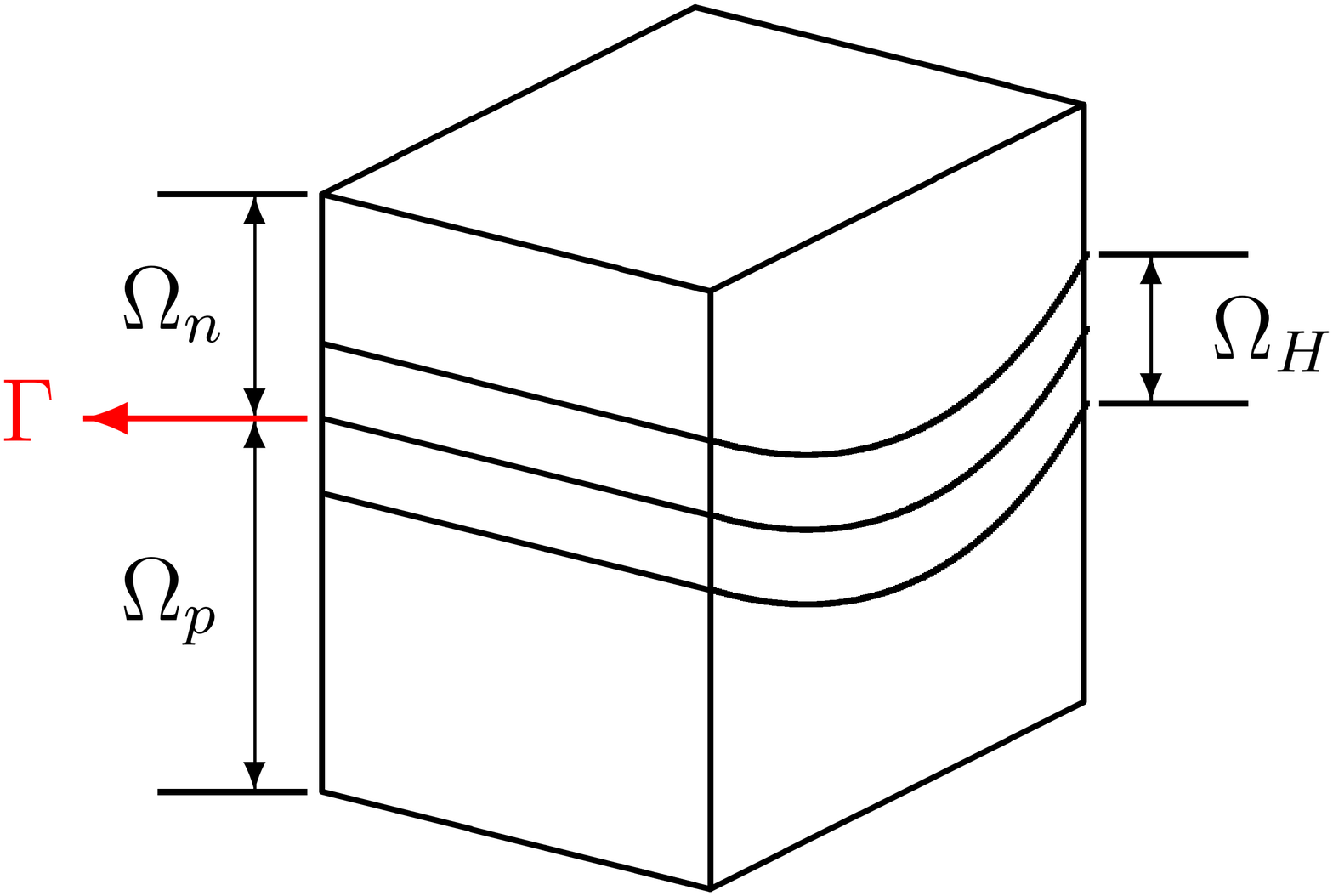}
\label{fig:OmegaH_3d}}
\subfigure[2D cross-section]
{\includegraphics[width=0.45\textwidth]{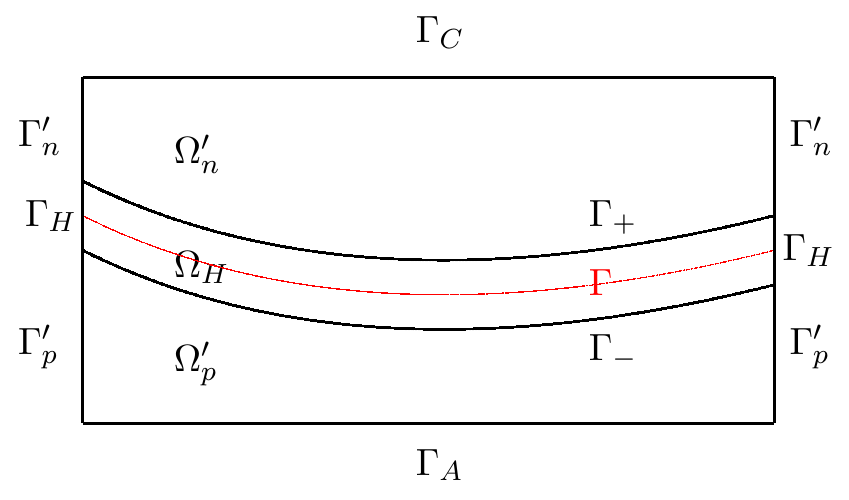}
\label{fig:OmegaH_2d}}
\caption{Geometry of the cell bulk and interface region.}
\label{fig:fig_OmegaH}
\end{figure}

Following~\cite{barker2003,williams_th,williams2008}, 
it is convenient, for modeling purposes, to 
associate with the interface $\Gamma$ the subregion 
$\Omega_H \subset \Omega$ depicted
in Fig.~\ref{fig:OmegaH_3d} and defined as follows.
For each point $\vect{x} \in \Gamma$, 
let $\vect{t}_{\vect{x}}=\left\{ \vect{x}+\xi \boldsymbol{\nu}_{\Gamma}
\left( \vect{x}\right) :\ \left\vert \xi \right\vert <H\right\}$ 
be the \lq\lq thickness\rq\rq{} associated with $\vect{x}$. Then, set
\begin{equation}\label{eq:omegaH}
\Omega _{H}=\bigcup\limits_{\vect{x} \in \Gamma }
\vect{t}_{\vect{x}}=\left\{ \vect{y}\in \Omega :
\,\mathrm{dist}(\vect{y},\Gamma)<H\right\}.
\end{equation}

The subregion $\Omega_H$ %(denoted \lq\lq physical interface\rq\rq) 
is thus a 3D thin layer of thickness $2 H$ surrounding $\Gamma$ 
which represents the device volumetric portion where the 
dissociation and recombination mechanisms of Sect.~\ref{sec:NBHJSCs} 
are assumed to occur. It is worth noting that 
the width $H$ is, in general, an unknown of the physical problem.
As such, it depends on $\vect{x}$ and $t$, it may assume different values 
in the two material phases of the photoactive
layer and  might in principle locally depend on the electric field.
According to the data provided 
in~\cite{barker2003,williams_th,williams2008}, we assume
for simplicity $H$ to be a constant quantity.
Based on the definition~\eqref{eq:omegaH}, we can introduce the two
portions $\Omega_{n}^{\prime}=\Omega_{n}\setminus \overline{\Omega}_{H}$ and
$\Omega_{p}^{\prime}=\Omega_{p}\setminus \overline{\Omega}_{H}$
%, in such a way that \meo{$\overline{\Omega} = \overline{\Omega}_{n}^{\prime }\cup %\overline{\Omega}_{H}\cup 
%\overline{\Omega}_{p}^{\prime }$ (CI VOGLIONO LE CHIUSURE?)}
(see Fig.~\ref{fig:OmegaH_2d}). 
Consistently, we also introduce the boundary portions $\Gamma_n^{\prime}$
and $\Gamma_p^{\prime}$ and set $\Gamma _{\pm }=\left\{ \vect{x} \pm H 
\boldsymbol{\nu}_{\Gamma}\left( \vect{x}\right) :\ \vect{x} \in \Gamma \right\}$,
in such a way that
$\partial \Omega _{n}^{\prime }=\Gamma _{C}\cup \Gamma _{+}\cup \Gamma
_{n}^{\prime }$,
$\partial \Omega _{p}^{\prime }=\Gamma _{A}\cup \Gamma _{-}\cup \Gamma
_{p}^{\prime }$ and $\partial \Omega _{H}=\Gamma _{+}\cup \Gamma _{-}\cup
\Gamma _{H}$, where $\Gamma_H = (\Gamma_n \cup \Gamma_p) 
\setminus (\Gamma_n^{\prime} \cup \Gamma_p^{\prime})$.
% By definition, we also have $\partial \Omega =
% \partial \Omega _{n}^{\prime } \cup \partial \Omega _{p}^{\prime }
% \cup \partial \Omega _{H}$.
Notice that, unlike $\Gamma$, the surfaces $\Gamma_{-}$ 
and $\Gamma_{+}$ can be regarded as \lq\lq{mathematical}\rq\rq{} 
interfaces.

\subsection{Modeling Assumptions}\label{sec:basic_assumptions}
Let us denote by $e$, $P$, $n$ and $p$ the 
\emph{volumetric densities} of (singlet) excitons, 
polaron pairs, electrons and holes, respectively,
and by $\mathbf{J}_e$, $\mathbf{J}_P$, $\mathbf{J}_n$ and 
$\mathbf{J}_p$ the associated particle fluxes.
Based on the physical working principles of an OSC 
illustrated in Sect.~\ref{sec:NBHJSCs}, on 
the heterogeneous geometrical decomposition of the device
introduced in Sect.~\ref{sec:geometry} and the
extensive numerical simulations reported in~\cite{buxton2007}, 
we make the following modeling and geometrical assumptions:
\begin{description}
\item[A.1] excitons can be generated at any position in
the cell, so that $e=e(\vect{x},t)$ is a nonnegative function over all
the cell domain $\Omega$;
\item[A.2] electrons (holes) are not able to penetrate the
donor (acceptor) material beyond the interface layer $\Omega_{H}$,
so that electrons (holes) are nonnegative functions
over $\Omega_n^{\prime} \cup \Omega_H$
($\Omega_p^{\prime} \cup \Omega_H$), and are
identically equal to zero in 
$\Omega_p^{\prime}$ ($\Omega_n^{\prime}$);
\item[A.3] polarons are trapped and immobile 
in the interface region $\Omega_H$, so that $P$ 
is a nonnegative function over $\Omega_H$ and 
identically equal to zero in $\Omega_n^{\prime} 
\cup \Omega_p^{\prime}$;
\item[A.4] the OSC is in the  \lq\lq off\rq\rq{} state at $t=0^-$
(that is, before illumination),
so that the initial condition for all the involved 
densities is $e(\vect{x},0)=
P(\vect{x},0)=n(\vect{x},0)=p(\vect{x},0)=0$ for
all $\vect{x} \in \Omega$;
\item[A.5] the geometry of the device is an infinite periodic
repetition of the computational domain of Fig.~\ref{fig:fig_Omega},
so that periodic boundary conditions
are enforced for all variables on the lateral boundary of $\Omega$.
\end{description}

\subsection{Microscale Model}\label{sec:microscale_model}
In this section, we illustrate the microscale model we advocate 
in this work to be a mathematical representation of the 
functioning of an OSC.

\begin{subequations}\label{eqs:excitons}
We take excitons to obey:
\begin{equation}\label{eq:micro_excitons}
\dfrac{\partial e}{\partial t}+\nabla \cdot \mathbf{J}_{e}
= S_e^B + S_e^H \qquad \text{in }\Omega
\end{equation}
where we use assumptions A.1 and A.3 to define
\begin{equation}\label{eq:micro_excitons_terms}
S_e^B = Q - \dfrac{e}{\tau_{e}}  \qquad \text{in }\Omega
\end{equation}
and
\begin{equation}\label{eq:micro_excitons_terms_interface}
S_e^H = \left\{
\begin{array}{ll}
0 & \qquad \text{in }\Omega_{n}^{\prime} \cup \Omega_p^\prime \\[1ex]
\eta k_{rec}P -\dfrac{e}{\tau _{diss}} & \qquad \text{in }\Omega _{H}.
\end{array}
\right.
\end{equation}
The superscripts $^B$ and $^H$ represent the fact that the corresponding 
volumetric production terms are defined in the bulk and
interface regions, respectively.
The term $Q$ denotes the rate at which 
excitons are generated by photon absorption and is henceforth assumed to be 
a nonnegative given function of time and position while
$\tau_{e}$ is the exciton lifetime in the bulk materials.
In the interface region $\Omega_H$ additional dissociation and recombination
mechanisms are taken into account and
$\tau_{diss}^{-1}$  and $\eta k_{rec}$ represent
the rate constants for the transition of excitons to the polaron state and
that of polarons back to the exciton state, respectively.
In particular $k_{rec}$ denotes the total rate of polaron
recombination events and $0 \leq \eta \leq 1$ the fraction of such
events which produce a singlet exciton.
As excitons have zero
net charge, their flux is driven by diffusion forces only, {\it i.e.}\
the flux density may be expressed as 
\begin{equation}\label{eq:e_fl}
\mathbf{J}_{e}=-D_{e}\nabla e \qquad \text{in }\Omega
\end{equation}
$D_{e}$ being the \emph{exciton diffusion coefficient}.
At the contacts we assume perfect exciton quenching~\cite{peumans2003}
so that 
\begin{equation}\label{eq:e_quenching}
e = 0 \qquad \text{on } \Gamma_{C} \cup \Gamma_{A}.
\end{equation}
%while homogeneous Neumann conditions are assumed 
%on $\Gamma_n \cup \Gamma_p$, 
\end{subequations}

\begin{subequations}\label{eqs:electronsandholes}
Because of assumption A.2, 
the following equations hold for electrons:
\begin{equation}\label{eq:n_b}
\left\{
\begin{array}{ll}
\dfrac{\partial n}{\partial t}+ %\dfrac{1}{q} 
\nabla \cdot \mathbf{J}_{n} = S_{n,p}^H &
\qquad \text{in } \Omega \setminus \Omega_{p}^{\prime} \\[1ex]
n\equiv 0 & \qquad \text{in }\Omega _{p}^{\prime } 
\end{array}
\right.
\end{equation}
and for holes:
\begin{equation}\label{eq:p_b}
\left\{
\begin{array}{ll}
\dfrac{\partial p}{\partial t}+ 
\nabla \cdot \mathbf{J}_{p} = S_{n,p}^H & 
\qquad \text{in } \Omega \setminus \Omega_{n}^{\prime} \\[1ex]
p\equiv 0 & \qquad \text{in }\Omega _{n}^{\prime} 
\end{array}
\right.
\end{equation}
where the term $S_{n,p}^H$ is defined as
\begin{equation}\label{eq:S_np^H}
S_{n,p}^H = \left\{
\begin{array}{ll}
k_{diss}P - \gamma n p & \qquad \text{in } \Omega_H\\[1ex]
0 & \qquad \text{in } \Omega _{n}^{\prime} \cup \Omega _{p}^{\prime}.
\end{array}
\right.
\end{equation}
Notice that $S_{n,p}^H$ is identically zero in the bulk region 
$\Omega _{n}^{\prime} \cup \Omega _{p}^{\prime}$ as 
electrons and holes can only recombine with each other,
so no recombination occurs where either of the two species is missing.
In the interface region $\Omega_H$ both electrons ad holes exist 
so the terms in $S_{n,p}^H$ take into account for polaron pair dissociation 
with $k_{diss}$ rate constant (see 
Sect.~\ref{sec:k_diss_model} for the model) and
bimolecular recombination with rate constant $\gamma$.
As electrons and holes each bear a non-zero net charge, 
their flux is driven by both
diffusion and electric drift forces~\cite{jerome96}, therefore:
\begin{equation}\label{eq:Jn}
\mathbf{J}_{n} = -D_{n}\nabla n - \mu _{n}n\mathbf{E}
\end{equation}
and
\begin{equation}\label{eq:Jp}
\mathbf{J}_{p} = 
-D_{p}\nabla p + \mu _{p}p\mathbf{E}
\end{equation}
where $\mathbf{E}$ is the electric field while
$D_{n}$, $\mu_{n}$ and $D_{p}$, $\mu_{p}$ are the 
diffusion coefficient
and mobility for electrons and holes, respectively.

Because of assumption A.2, the following boundary
conditions hold at the artificial interfaces separating the
donor and acceptor bulk phases from the thin slab region $\Omega_H$:
\begin{equation}\label{eq:n_noflux}
\boldsymbol{\nu}_{\Gamma} \cdot \mathbf{J}_{n} = 0
\qquad \mbox{on } \Gamma_{-}
\end{equation}
and
\begin{equation}\label{eq:p_noflux}
\boldsymbol{\nu}_{\Gamma} \cdot \mathbf{J}_{p} = 0
\qquad \mbox{on } \Gamma_{+}.
\end{equation}
At the contacts we impose the same Robin-type boundary conditions as
described in~\cite{scott_malliaras99,defalco}:
\begin{equation}\label{eq:n_robin}
- \kappa_n \boldsymbol{\nu} \cdot \mathbf{J}_{n} +
\alpha_n n = \beta_n \qquad \mbox{on } \Gamma_{C}
\end{equation}
and
\begin{equation}\label{eq:p_robin}
- \kappa_p \boldsymbol{\nu} \cdot \mathbf{J}_{p} +
\alpha_p p = \beta_p  \qquad \mbox{on } \Gamma_{A},
\end{equation}
where $\kappa_n$, $\kappa_p$, $\alpha_n$, $\alpha_p$
$\beta_n$, $\beta_p$ are nonnegative coefficients.
%while no-flux boundary conditions are taken 
%to hold on $\Gamma_n \cup \Gamma_p$, this expressing
%again the fact that all photogenerated carriers remain 
%confined within the device.
\end{subequations}

\begin{subequations}\label{eqs:poisson}
The \emph{electric field} $\mathbf{E}$ in~\eqref{eq:Jn} and~\eqref{eq:Jp} 
is connected to the \emph{electric potential} $\varphi$ by the
quasi-static approximation
\begin{equation}\label{eq:efield}
\mathbf{E}=-\nabla \varphi  \qquad \text{in }\Omega
\end{equation}
and satisfies the Poisson equation
\begin{equation}\label{eq:Poisson}
\nabla \cdot \left( \varepsilon \mathbf{E} \right) = \rho \qquad
\mbox{in } \Omega
\end{equation}
where $\rho$ is the space charge density in the device.
Using assumption A.2, the piecewise smooth definition of $\rho$ 
turns out to be:
\begin{equation}\label{eq:spacechargedensity}
\rho = q (p - n) =
\left\{
\begin{array}{ll}
-q n    &  \quad \mbox{in } \Omega_n^{\prime} \\
q (p-n) & \quad \mbox{in }  \Omega_H \\ 
+q p    & \quad \mbox{in } \Omega_p^{\prime},
\end{array}
\right.
\end{equation}
% -q n \chi_{} + q p \chi_{\Omega_p^{\prime}}
% + q (p -n) \chi_{\Omega_H}.
% \end{equation}
% In~\eqref{eq:spacechargedensity}, 
% we have denoted, for any set $\mathcal{S} \subseteq \Omega$, 
% by $\chi_{\mathcal{S}}$ the characteristic function of $\mathcal{S}$, such
% that $\chi_{\mathcal{S}}(\vect{x})=1$ if $\vect{x} \in \mathcal{S}$
% and $\chi_{\mathcal{S}}(\vect{x})=0$ otherwise.
% Notice that the piecewise smooth definition of $\rho$
% is consistent with  on the carrier densities.}

$q$ denoting the quantum of charge.
The \emph{electric permittivity} $\varepsilon$ is equal to $\varepsilon_r \varepsilon_0$, 
$\varepsilon_r$ and $\varepsilon_0$ being the relative material 
and vacuum permittivities, respectively, with 
$\varepsilon_{r} = \varepsilon_{r,a}$
in the acceptor phase and $\varepsilon_{r} = \varepsilon_{r,d}$ 
in the donor phase,
so that $\varepsilon$ may be discontinuous across the interface $\Gamma$.
Dirichlet boundary conditions for the electric potential are set 
at the contacts $\Gamma_A$ and $\Gamma_C$, as follows
\begin{equation}\label{eq:bc_phi_gamma_c}
\varphi = 0
\qquad \mbox{on } \Gamma_{C}
\end{equation}
and
\begin{equation}\label{eq:bc_phi_gamma_a}
\varphi =  V_{appl} + V_{bi}
\qquad \mbox{on } \Gamma_{A}
\end{equation}
where $V_{bi} = (\Phi_A-\Phi_C)/q$ is the built-in voltage
of the cell, $\Phi_A$ and $\Phi_C$ are 
the contact metal work functions while 
$V_{appl}$ is the externally applied voltage.
%On $\Gamma_n \cup \Gamma_p$ 
%homogeneous Neumann boundary conditions are assumed to hold,
%this expressing the fact that the electric field force lines
%start and close at the contacts.
\end{subequations}

\begin{subequations}\label{eq:polarons}
Because of assumption A.3, the flux $\mathbf{J}_P$ is identically
equal to zero in all $\Omega$ and for all $t\geq 0$, and
polarons satisfy the following ODE in the interface region 
\begin{equation}\label{eq:polaronkinetics}
\dfrac{\partial P}{\partial t}=\dfrac{e}{\tau _{diss}}+\gamma np-
\left(k_{diss}+k_{rec}\right) P  \qquad \text{in }\Omega _{H} 
\end{equation}
while their density is identically zero in the bulk
\begin{equation}\label{eq:polaronkinetics2}
P\equiv 0  \qquad \text{in }\Omega _{n}^{\prime }\cup 
\Omega _{p}^{\prime}. 
\end{equation}
\end{subequations}

\subsection{Micro-to-Macro Scale Transition}
\label{sec:macroscale_model}

The microscale model of a bilayer OSC described 
in Sect.~\ref{sec:microscale_model} can be subdivided into
three distinct groups of equations:
\begin{description}
\item[1)] parabolic PDEs enforcing mass conservation of excitons, 
electrons and holes; 
\item[2)] an ODE
describing the kinetics of photogenerated polaron pairs;
\item[3)] an elliptic constraint enforcing Gauss theorem 
in differential form to be satisfied at each
time $t>0$ throughout the whole cell domain.
\end{description}

The markedly spatially heterogeneous nature of the problem 
may be quite impractical for numerical  
simulation, in particular when devices with complex interface
morphology in multiple spatial dimensions are considered. 
For this reason, in this section we propose a 
scale transition procedure which allows us to derive 
a macroscale model that is more amenable to numerical treatment. 
Other examples of multiscale mathematical approaches that
are based on the scale separation concept and scale transition
can be found in~\cite{Jaffre2005,transitions,Mori2009}.

To construct our multiscale model of an OSC, we abandon the 
perspective focused at the nanoscopic characteristic level 
adopted so far, and prefer to look at 
the cell from a \lq\lq larger\rq\rq{} distance.
By doing so, necessarily, we loose control of the details
(i.e, we cannot distinguish the region $\Omega_H$ from the
two bulk regions $\Omega_n$ and $\Omega_p$), but, at the same time,
we gain the advantage of not needing to resolve the interfacial
bulk region across $\Gamma$. The resulting macroscale problem 
is thus posed in the {\em partitioned} domain 
$\Omega \setminus \Gamma \equiv \Omega_n \cup \Omega_p$ (as a matter
of fact, we are still able to neatly distinguish the interface
separating the two material phases!) {\em without including} 
the interfacial production terms $S^H_{(\cdot)}$ in the mass balance 
and kinetics equations 1) and 2) introduced above. 

Of course, we cannot simply limit ourselves to 
neglecting these latter terms, rather, we do need to incorporate 
their effects, in the macroscale model, in an alternative way. 
For this, the simplest approach to micro-to-macro scale 
transition consists of replacing $S^H_{(\cdot)}$, at each 
point of $\Gamma$ and for each time level, with its {\em average} 
$\sigma^H_{(\cdot)}$ across the thickness of 
$\Omega_H$ in the normal direction, 
and then, of using $\sigma^H_{(\cdot)}$ as a source term for suitable
flux transmission conditions, to be enforced on the interface $\Gamma$ 
in the case of mass balance equations. In the case of 
polaron pair equation, the averaging procedure automatically 
transforms the volumetric kinetics balance within $\Omega_H$
into a {\em surface} kinetics balance over $\Gamma$.
In any case, the scale transition results in the introduction
of suitable interfacial terms that replace 
in a \lq\lq lumped\rq\rq{} manner the volumetric 
dissociation/generation phenomena microscopically occurring 
in $\Omega_H$. Having characterized the averaging procedure 
for equations 1) and 2), the (macroscale) differential Gauss
theorem 3) remains automatically (formally) unchanged and is 
expressed in terms of the (macroscale) space charge density
as in Eqns.~\eqref{eqs:poisson}.

\subsubsection{Derivation of the Macroscale Equations}
\label{sec:derivation_of_macro_model}

\begin{subequations}\label{eqs:excitons_macro}
The macroscale model for excitons reads:
\begin{equation}\label{eq:macro_excitons}
\dfrac{\partial e}{\partial t}+\nabla \cdot \mathbf{J}_{e}
= S_e^B  \qquad \text{in }\Omega \setminus \Gamma
\end{equation}
with
\begin{equation}\label{eq:macro_excitons_terms}
S_e^B = Q - \dfrac{e}{\tau_{e}}  \qquad \text{in }\Omega
\setminus \Gamma
\end{equation}
and subject to the interface conditions:
\begin{equation}\label{eq:transmission_excitons}
\left\{
\begin{array}{ll}
\jmp{\boldsymbol{\nu}_\Gamma \cdot \mathbf{J}_{e}} = 
\sigma_{e}^H & \qquad \text{on } \Gamma \\[1ex]
\jmp{e} = 0 & \qquad \text{on } \Gamma
\end{array}
\right.
\end{equation}
where $\jmp{f}:=f_n - f_p$ denotes for any function 
$f:\Omega \rightarrow \mathbb{R}$ the jump of $f$ across the
interface $\Gamma$, $f_n$ and $f_p$ being the traces on $\Gamma$ 
of the restrictions of $f$ from $\Omega_n$ and $\Omega_p$, respectively.
The continuity of $e$ at the interface is a requirement 
consistent with the elliptic regularity of both microscale
and macroscale problems.

The interfacial source term $\sigma_e^H$ is defined as
\begin{equation}\label{eq:sigmaH_e}
\sigma_e^H =
\int_{-H}^{H} \left(\eta k_{rec} P - \dfrac{e}{\tau_{diss}}\right) \, d\xi =
\eta k_{rec} \int_{-H}^{H} P \, d\xi - 
\dfrac{1}{\tau_{diss}}\int_{-H}^{H} e \, d\xi \simeq 
\eta k_{rec} \widetilde{P} - \dfrac{2H}{\tau_{diss}}\,e\vert_{\Gamma}.
\end{equation}
In the above relation, $e\vert_{\Gamma}$ is the (single--valued)
trace of $e$ over $\Gamma$, while
$\widetilde{P}$ is the {\em areal density} of the bonded pairs, defined as
\begin{equation}\label{eq:polaron_areal_density}
\widetilde{P}(\vect{x},t) =
\int_{-H}^{H} P(\vect{x} + \xi \vect{t}_{\vect{x}},t) \, d\xi 
\qquad \forall \vect{x} \in \Gamma
\end{equation}
and the midpoint quadrature rule is used for approximating 
the third integral in~\eqref{eq:sigmaH_e}.
The macroscale model for excitons is completed by the
constitutive relation~\eqref{eq:e_fl} for exciton flux density
and by the perfect exciton quenching 
boundary conditions~\eqref{eq:e_quenching}.
\end{subequations}

\begin{subequations}\label{eqs:electrons_macro}
The macroscale model for electrons reads:
\begin{equation}\label{eq:macro_electrons}
\left\{
\begin{array}{ll}
\dfrac{\partial n}{\partial t}+\nabla \cdot \mathbf{J}_{n}
= 0  & \qquad \text{in }\Omega_n \\[1ex]
n \equiv 0 & \qquad \text{in }\Omega_p
\end{array} 
\right.
\end{equation}
subject to the interface/boundary condition
\begin{equation}\label{eq:transmission_electrons}
\boldsymbol{\nu}_\Gamma \cdot \mathbf{J}_{n} = 
\sigma_{n,p}^H \qquad \text{on } \Gamma.
\end{equation}
The interfacial source term $\sigma_{n,p}^H$ is defined as
\begin{equation}\label{eq:sigmaH_n}
\sigma_{n,p}^H = \int_{-H}^{H}
\left(k_{diss} P - \gamma n p \right) \, d \xi \simeq
k_{diss}\vert_\Gamma \int_{-H}^{H} P\, d \xi 
- \int_{-H}^{H} \gamma n p \, d \xi \simeq 
k_{diss}\vert_\Gamma \widetilde{P} - 2H\, \gamma\vert_\Gamma\, 
n\vert_\Gamma\, p\vert_\Gamma 
\end{equation}
where definition~\eqref{eq:polaron_areal_density} 
is used in the first integral while 
the midpoint quadrature rule is again used to approximate 
the third integral in~\eqref{eq:sigmaH_n}.
The macroscale model for electrons is completed by the
constitutive relation~\eqref{eq:Jn} for electron flux density
and by the Robin-type boundary condition~\eqref{eq:n_robin}.
\end{subequations}

\begin{subequations}\label{eqs:holes_macro}
Proceeding in a completely analogous manner as done with
electrons, the macroscale model for holes reads:
\begin{equation}\label{eq:macro_holes}
\left\{
\begin{array}{ll}
\dfrac{\partial p}{\partial t}+\nabla \cdot \mathbf{J}_{p}
= 0  & \qquad \text{in }\Omega_p \\[1ex]
p \equiv 0 & \qquad \text{in }\Omega_n
\end{array} 
\right.
\end{equation}
subject to the interface/boundary condition
\begin{equation}\label{eq:transmission_holes}
\boldsymbol{\nu}_\Gamma \cdot \mathbf{J}_{p} = 
-\sigma_{n,p}^H \qquad \text{on } \Gamma.
\end{equation}
The macroscale model for holes is completed by the
constitutive relation~\eqref{eq:Jp} for hole flux density
and by the Robin-type boundary condition~\eqref{eq:p_robin}.

The conditions~\eqref{eq:transmission_electrons} and~~\eqref{eq:transmission_holes} assume an interesting
physical meaning upon introducing the electron and hole 
{\em current densities}, defined respectively 
as $\vect{j}_n:=-q \mathbf{J}_{n}$ and 
$\vect{j}_p:=+q \mathbf{J}_{p}$, and the 
{\em total (conduction) current density}
$\vect{j}:=\vect{j}_{n} + \vect{j}_{p}$. Recalling that
$n=0$ ($p=0$) in $\Omega_p$ ($\Omega_n$), we have:
$$
\vect{j} = 
\left\{
\begin{array}{ll}
\vect{j}_{n} & \qquad \mbox{in } \Omega_n \\[1ex]
\vect{j}_{p} & \qquad \mbox{in } \Omega_p
\end{array}
\right.
$$
from which we get
\begin{equation}\label{eq:current_conservation}
\jmp{\boldsymbol{\nu}_\Gamma \cdot \vect{j}} = 0 
\qquad \mbox{on } \Gamma,
\end{equation}
that expresses the property of {\em current conservation} 
across the interface $\Gamma$.
% \begin{equation}\label{eq:sigmaH_p}
% \sigma_p^H = \int_{-H}^{H}
% \left(k_{diss} P - \gamma n p \right) \, d \xi \simeq
% k_{diss}\vert_\Gamma \int_{-H}^{H} P\, d \xi 
% - \int_{-H}^{H} \gamma n p \, d \xi \simeq 
% k_{diss}\vert_\Gamma \widetilde{P} - 2H\, \gamma\vert_\Gamma\, 
% n\vert_\Gamma\, p\vert_\Gamma.
% \end{equation}
\end{subequations}

\begin{subequations}\label{eqs:polarons_macro}
Integration of~\eqref{eq:polaronkinetics} across the interface thickness
yields the following macroscale model for the areal density 
of polaron pairs
\begin{equation}\label{eq:macro_polarons}
\dfrac{\partial \widetilde{P}}{\partial t} = \sigma_P^H 
\qquad \text{on }\Gamma
\end{equation}
where
\begin{equation}\label{eq:sigmaH_P}
\sigma_P^H = 
\dfrac{2H}{\tau _{diss}} e\vert_\Gamma 
+ 2H\,\gamma\vert_\Gamma\, n\vert_\Gamma\, p\vert_\Gamma
- \left(k_{diss}\vert_\Gamma+k_{rec}\right) \widetilde{P}.
\end{equation}
% consistently with the conservation of mass principle within the 
% interfacial region $\Omega_H$.
\end{subequations}

\begin{subequations}\label{eqs:poisson_macro}
The macroscale model for the differential Gauss
theorem is expressed by the following Poisson problem in
heterogeneous form:
\begin{equation}\label{eq:macro_potential}
\nabla \cdot \left( \varepsilon \mathbf{E} \right) = \rho \qquad
\mbox{in } \Omega \setminus \Gamma
\end{equation}
with
\begin{equation}\label{eq:macro_rho}
\rho = 
\left\{
\begin{array}{ll}
- q n & \quad \text{in } \Omega_n \\[1ex]
+ q p & \quad \text{in } \Omega_p 
\end{array}
\right.
\end{equation}
and subject to the interface conditions:
\begin{equation}\label{eq:transmission_poisson}
\left\{
\begin{array}{ll}
\jmp{\boldsymbol{\nu}_\Gamma \cdot \varepsilon \mathbf{E}} = 
0 & \qquad \text{on } \Gamma \\[1ex]
\jmp{\varphi} = 0 & \qquad \text{on } \Gamma.
\end{array}
\right.
\end{equation}
Two remarks are in order with system~\eqref{eqs:poisson_macro}. First, 
we notice that the Gauss theorem in differential 
form~\eqref{eq:macro_potential} looks formally identical
to the corresponding microscale formulation~\eqref{eq:Poisson}, 
the difference between the two methodologies being in 
the definition of the space charge density $\rho$
(compare~\eqref{eq:spacechargedensity} with~\eqref{eq:macro_rho}).
Second, the transmission conditions~\eqref{eq:transmission_poisson} 
express the physical fact that the normal component of
the electric displacement vector and the electric potential
do not experience any discontinuity at the material interface, as 
is the case of the microscale formulation.
\end{subequations}

\subsubsection{Summary of the Macroscale Model}
\label{sec:complete_microscale_model}
For sake of convenience, we summarize below 
the macroscale model of an OSC written in primal form:
\begin{subequations}\label{eq:lumped_model}
\begin{align}
\label{e_red}
&\left\{
\begin{array}{ll}
\dfrac{\partial e}{\partial t}-\nabla \cdot \left( D_{e}\nabla e\right) =
Q-\dfrac{e}{\tau _{e}} & \qquad \text{in }\Omega _{n}\cup \Omega _{p} 
\equiv \Omega \setminus \Gamma\\[1ex]
\jmp{e} = 0, \;\;
\jmp{-\boldsymbol{\nu}_{\Gamma} \cdot D_{e}\nabla e } =\eta
k_{rec} \widetilde{P}-\dfrac{2H}{\tau _{diss}}e%
& \qquad \text{on }\Gamma, \\[1ex]
e = 0 & \qquad \mbox{on } \Gamma_{C} \cup \Gamma_A, \\[0.5ex]
%\mathbf{J}_{e} \cdot \boldsymbol{\nu} = 0 
%& \qquad \mbox{on } \Gamma_{n}, \\
e(\vect{x},0) = 0, \;\;
& \qquad \forall \vect{x} \in \Omega,
\end{array}%
\right.\\[3mm]
\label{Polaron_red}
&\left\{
\begin{array}{ll}
\dfrac{\partial \widetilde{P}}{\partial t}=\dfrac{2H}{\tau _{diss}}%
e+2H\gamma np-\left( k_{diss}+k_{rec}\right) \widetilde{P} & 
\qquad \text{on } \Gamma,\\[2ex]
\widetilde{P}(\vect{x},0) = 0, \;\;
& \qquad \forall \vect{x} \in \Gamma,
\end{array}%
\right.\\[3mm]
\label{n_red}
&\left\{
\begin{array}{ll}
\dfrac{\partial n}{\partial t}- 
\nabla \cdot \left(D_{n}\nabla n - \mu _{n}n\nabla \varphi\right) = 0 & 
\qquad \text{in }\Omega _{n} \\[1.5ex]
\boldsymbol{\nu}_{\Gamma} \cdot \left(D_{n}\nabla n - \mu _{n}n\nabla \varphi\right) = 
\mathbf{-}k_{diss}\widetilde{P}+2H\gamma np & \qquad \text{on }\Gamma, \\[0.5ex]
\kappa_n \boldsymbol{\nu} \cdot \left(D_{n}\nabla n - 
\mu _{n}n\nabla \varphi \right) +
\alpha_n n = \beta_n & \qquad \mbox{on } \Gamma_{C}, \\[0.5ex]
%\mathbf{J}_{n} \cdot \boldsymbol{\nu} = 0 & 
%\qquad \mbox{on } \Gamma_{n}, \\
n(\vect{x},0) = 0, \;\;
& \qquad \forall \vect{x} \in \Omega,
\end{array}
\right.\\[3mm]
\label{p_red}
&\left\{
\begin{array}{ll}
\dfrac{\partial p}{\partial t}-
\nabla \cdot \left(D_{p}\nabla p + \mu _{p} p\nabla \varphi\right) =0 & 
\qquad \text{in }\Omega _{p} \\[1.5ex]
- \boldsymbol{\nu}_{\Gamma} \cdot \left(D_{p}\nabla p + 
\mu _{p} p\nabla \varphi\right) = 
\mathbf{-}k_{diss}\widetilde{P}+2H\gamma np & \qquad \text{on }\Gamma, \\[0.5ex]
\kappa_p \boldsymbol{\nu} \cdot 
\left(D_{p}\nabla p + \mu _{p} p\nabla \varphi\right) +
\alpha_p p = \beta_p & \qquad \mbox{on } \Gamma_{A}, \\[0.5ex]
%\mathbf{J}_{p} \cdot \boldsymbol{\nu} = 0 & 
%\qquad \mbox{on } \Gamma_{p}, \\
p(\vect{x},0) = 0, \;\;
& \qquad \forall \vect{x} \in \Omega,
\end{array}%
\right.\\[3mm]
\label{Poisson_red}
&\left\{
\begin{array}{ll}
-\nabla \cdot \left( \varepsilon \nabla \varphi \right) =-q \, n & 
\qquad \text{in } \Omega_n  \\[0.5ex]
-\nabla \cdot \left( \varepsilon \nabla \varphi \right) =+q \, p & 
\qquad \text{in } \Omega_p  \\[0.5ex] 
\jmp{\varphi} = \jmp{\boldsymbol{\nu}_{\Gamma} \cdot \varepsilon \nabla \varphi } 
=0 & \qquad \text{on }\Gamma, \\[0.5ex]
\varphi = 0 & \qquad \mbox{on } \Gamma_{C}, \\[0.5ex]
\varphi = V_{appl} + V_{bi}
& \qquad \mbox{on } \Gamma_{A}.
%-\varepsilon \nabla \varphi \cdot \boldsymbol{\nu} = 0 & 
%\qquad \mbox{on } \Gamma_{n} \cup \Gamma_p.
\end{array}
\right.
\end{align}
\end{subequations}
System~\eqref{eq:lumped_model} is completed by periodic 
boundary conditions on $\Gamma_n \cup \Gamma_p$, as stated 
in assumption A.5.
For the physical models of the coefficients 
in system~\eqref{eq:lumped_model} we refer 
to~\cite{barker2003,gill1972,horowitz1998}, except for the
description of the polaron dissociation rate constant $k_{diss}$
which is addressed in detail in Sect.~\ref{sec:k_diss_model}.
In particular, for the carrier mobilities, we neglect the 
effect of energetic disorder, so that they can be
assumed to depend only on the electric field magnitude, 
according to the Poole-Frenkel model.
As for diffusivities, 
in the computations of Sect.~\ref{sec:num_res}, Einstein relations
\begin{equation}\label{eq:einsteinrelation}
D_n = (K_B T/q) \mu_n, \qquad 
D_p = (K_B T/q) \mu_p
\end{equation}
are assumed to hold, although the proposed multiscale formulation 
remains unchanged if such an assumption is removed.
In~\eqref{eq:einsteinrelation}, 
$K_B$ is Boltzmann's constant and $T$ is the absolute temperature.
Finally, for the bimolecular recombination rate constant $\gamma$
a Langevin-type relation is used~\cite{barker2003}.

\subsection{Model for the Polaron Dissociation Rate}
\label{sec:k_diss_model}

Numerical simulations as those reported in Sect.~\ref{sec:num_res}
show  that the polaron dissociation rate $k_{diss}$ has a significant
impact on the cell photoconversion efficiency, for this reason we
devote  this entire section to modeling the dependence of $k_{diss}$ on the
electric field and on the morphology of the material interface.  
A commonly used polaron dissociation rate model is the
Braun-Onsager model~\cite{Braun1984} which is derived assuming the OSC bulk 
to be a homogeneous medium and takes into account only the magnitude of the electric field. 
In~\cite{barker2003} the authors propose a model for
$k_{diss}(\mathbf{E})$ tailored for bilayer devices which is derived
by performing an average over the admissible range of the escape angle 
relative to the electric field direction. In this latter model, the
electric field is assumed to be always directed orthogonally to the
interface consistently with the planar geometry of the device considered therein.
The authors of~\cite{williams2008} apply the dissociation rate model 
of~\cite{barker2003} to more complex
geometries by performing an  average along the interface of the field
component normal to the contacts which amounts to neglecting  the
effect of local electric field orientation. 

To construct a novel model which also takes into account this latter
effect we repeat the derivation of~\cite{barker2003} with two
differences. The first difference is that of removing the assumption that 
the field is normal to the interface. 
The second difference is that of considering a limited range of
admissible escape directions to account for the physical fact that
polaron pairs tend to be aligned with the gradient of the electron 
affinity due to the different materials in the two device subregions.

\begin{figure}[hbt!]
\begin{center}
\includegraphics[angle=0,width=.6\textwidth]{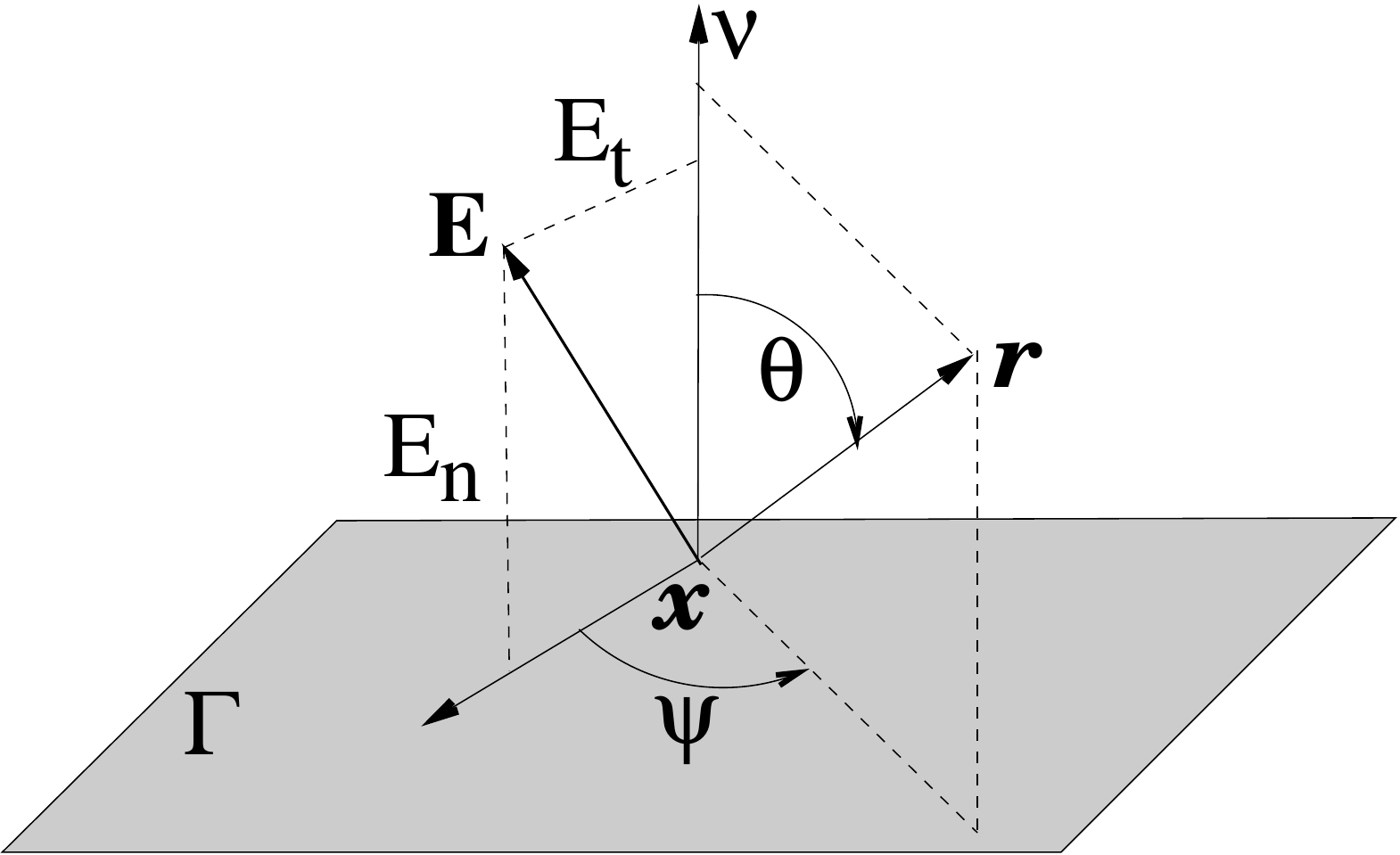}
% \qquad
% \includegraphics[width=.32\textwidth]{interf2.png}
\caption{Geometrical notation of the quantities involved in 
polaron dissociation at the material interface.}
\label{fig:interfaccia}
\end{center}
\end{figure}

Referring to Fig.~\ref{fig:interfaccia} for the geometrical notation,
we let  
\begin{equation}\label{eq:kdissmodel1}
k_{diss} (\mathbf{E}) = \displaystyle k_{diss} (0) \int_{0}^{2\pi}
d\psi \int_{0}^{\pi/2} w(\theta, \psi)\ \beta \left(
  \mathbf{E}\cdot\vect{r}\right)  d\theta, %\sin \left(\theta\right) 
\end{equation}
where $k_{diss}(0)$ is the zero-field dissociation rate
constant, $\vect{r}$ is the escape direction of the electron part 
of the polaron at the point $\vect{x} \in \Gamma$, 
$w$ is a nonnegative weight representing the probability
  distribution of admissible escape directions, and such that
$\int_{0}^{2\pi} d\psi  
\int_{0}^{\pi/2} w(\theta, \psi)\, d\theta = 1$, and
$\beta$ is an enhancement/suppression factor given by the 
Poole-Frenkel formula:
\begin{equation}\label{eq:enhancement}
\beta(z) = 
\left\{
\begin{array}{ll}
e^{-Az} & \qquad z \geq 0\\[1ex]
e^{2\sqrt{-Az}} & \qquad  z < 0,
\end{array}
\right.
\end{equation}
having set $A = (4\pi\varepsilon)^{-1} q^3 (K_b \,T)^{-2}$. 
The product $\mathbf{E} \cdot \vect{r}$ can be expressed in terms of
the normal component $E_{n}$ and the tangential component $E_{t}$ of
the electric field as 
\[
\mathbf{E} \cdot \vect{r} = E_{n} \cos \theta + E_{t}
\sin \theta \cos \psi.
\]
To specify an expression for $w$ we assume an escape 
direction $\vect{r}$ to be admissible
only when the angle it forms with respect to the normal unit vector $\boldsymbol{\nu}$ is not too large.
Indicating by $\theta_{max}$ the maximum admissible value for 
$\theta$ and allowing all 
admissible values to be equally likely, we obtain:
% 
% Upon assuming that $\vect{v}$ forms a maximum angle $\theta_{max}$
% with $\boldsymbol{\nu}$, the weight $w$ can be expressed as: 
$$
w(\theta, \psi) =
\left\{
\begin{array}{ll}
\displaystyle \frac{\sin \theta}{2 \pi (1-\cos \theta_{max})} 
& \qquad 0 < \theta \leq \theta_{max} \\[4mm]
0 & \qquad \theta_{max} < \theta \leq \displaystyle \frac{\pi}{2}.
\end{array}
\right.
$$
Two limits are of particular interest, 
$\theta_{max} \rightarrow 0^+$ and 
$\theta_{max}=\pi/2$. 

%%  the discontinuity
%% in electron affinity across the interface tends to align 
%% the principal axis of the geminate pair dipole to the local
%% interface normal vector. The reason why we do not assume 
%% {\em a priori} $\theta_{max} = 0$ is that we want to allow 
%% our model to account for statistical deviations from the expected value  
%% due to,  e.g., interface surface roughness and/or thermal vibration mechanisms.

%% The case $\theta_{max}=\pi/2$ is significant because it 
%% corresponds to the model proposed in Ref. [3]. 
%% Notice that taking $\theta_{max}=\pi/2$ amounts to 
%% completely neglecting the alignment behavior mentioned above,
%% and, thus, it is no surprise that such a choice might
%% lead to the overestimation of the dissociation rate
%% which is actually observed in Fig. 8(a).

In the first case, Eq.~\eqref{eq:kdissmodel1} can be checked to yield
\begin{equation}\label{eq:kdissmodel3}
k_{diss} (\mathbf{E}) =  k_{diss} (0)\ \beta \left(E_{n}\right).
\end{equation}
This corresponds to assuming that all geminate pairs are exactly aligned
with the interface normal unit vector, thus neglecting any possible variability 
in their orientation due to, {\it e.g.}\ interface surface roughness and/or thermal vibrations.
  
In the second case, Eq.~\eqref{eq:kdissmodel1} becomes
\begin{equation}\label{eq:kdissmodel2}
k_{diss} (\mathbf{E}) =  
\displaystyle k_{diss} (0) \int_{0}^{2\pi} d\psi 
\int_{0}^{\pi/2} \frac{\sin\theta}{2 \pi} \ 
\beta \left( \mathbf{E}\cdot\vect{r}\right) d\theta, %\sin \left(\theta\right)
\end{equation}
which, in the special case where $E_t=0$, coincides with 
Eqs.~(17)-(21) of~\cite{barker2003}.
Notice that if $E_t \neq 0$ 
the choice $\theta_{max} = \pi/2$ may overestimate the effective dissociation rate as it corresponds to 
completely neglecting the alignement of the geminate pairs with the electron affinity gradient.
This is observed to give rise to non-physical effects 
as shown by the simulations of Sect.~\ref{sec:comparison}.
Therefore, for practical purposes, the quantity $\theta_{max}$ 
should be used as a fitting parameter to be calibrated on experimental data.

\begin{figure}[!hbt]
	\centering
	\subfigure[Model~\eqref{eq:kdissmodel3}]{
		\includegraphics[width=.475\textwidth]{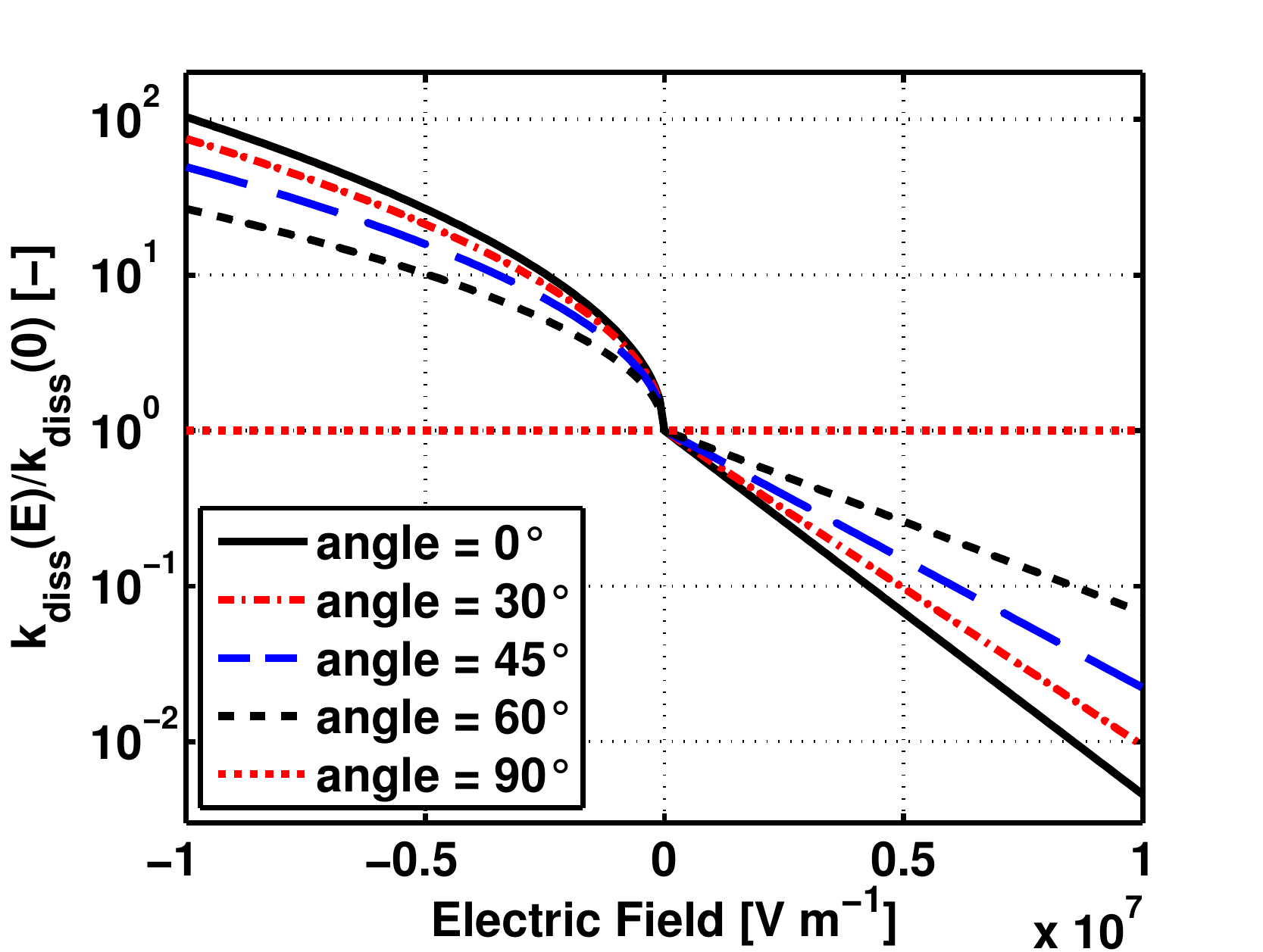}
	}
	\subfigure[Model~\eqref{eq:kdissmodel2}]{
		\includegraphics[width=.475\textwidth]{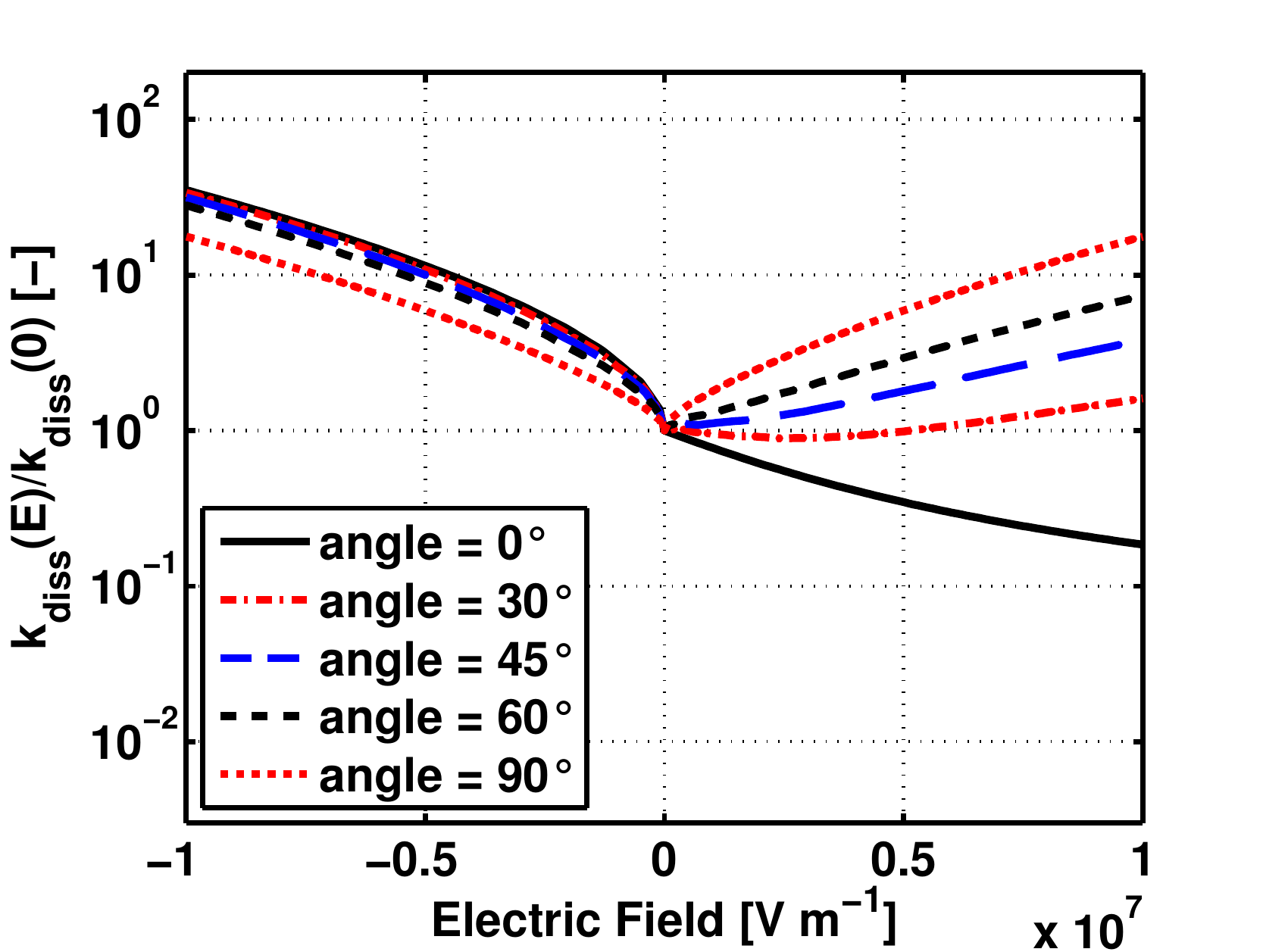}
	}
	\caption{Comparison between models~\eqref{eq:kdissmodel3}
          and~\eqref{eq:kdissmodel2} 
	for various angles between $\mathbf{E}$ and
        $\boldsymbol{\nu}$.} 
	\label{fig:kdissvsE}
\end{figure}

Fig.~\ref{fig:kdissvsE} shows the dissociation
rate constant (normalized to $k_{diss}(0)$)
computed by model~\eqref{eq:kdissmodel3} (left) 
and~\eqref{eq:kdissmodel2} (right) for several values of the angle
between $\mathbf{E}$ and $\boldsymbol{\nu}$ and having set $T= 300 \, 
\kelvin$ 
and $\varepsilon_r=4$.
We notice that the dissociation rate computed by
model~\eqref{eq:kdissmodel2}  
has a significantly smaller range of variability than predicted by
model~\eqref{eq:kdissmodel3}.
A possible explanation to this difference is related to the
smoothing operated by the integral in~\eqref{eq:kdissmodel2}.
The higher variability of the dissociation rate translates into an higher  
sensitivity of model~\eqref{eq:kdissmodel3} to the inclination of the electric field with 
respect to the interface normal as will be further discussed in the numerical
results section when commenting Fig.~\ref{fig:angles_right}.
%% If the electric field lies  
%% in the tangent plane at point $\vect{x}$ of $\Gamma$
%% (angle $=90^{\circ}$, $E_n=0$), the resulting curve from
%% model~\eqref{eq:kdissmodel2} is an even function of $E_t$,  
%% similarly to what predicted by classical Onsager's dissociation 
%% theory~\cite{onsager,Braun1984}, 
%% while in the case of model~\eqref{eq:kdissmodel3} 
%% the normalized dissociation does not depend on the electric field
%% because $\beta(E_n)=1$.
A discussion of the impact of~\eqref{eq:kdissmodel3} 
and~\eqref{eq:kdissmodel2} on the model predictions 
will be carried out in Sect.~\ref{sec:comparison}.

%%% Local Variables: 
%%% mode: latex
%%% TeX-master: "paper_osc_multiscale"
%%% End: 

\section{Numerical Approximation}
\label{sec:numer_approx}
In this section we describe the numerical techniques used to solve the
mathematical models introduced in Sects.~\ref{sec:microscale_model}
and \ref{sec:macroscale_model}. The full details of the 
discrete system of linear algebraic equations resulting from problem 
approximation are postponed to~\ref{sec:appendix}.
As for the simulation of the model in the transient regime
carried out in Sect.~\ref{sec:comparison_macro_micro},
we have adapted to the case at hand
the numerical method described in~\cite{defalco} based on Rothe's
method and on  
the use of adaptive Backward Differentiation Formulas (BDF).
In the steady-state simulations illustrated in
Sects.~\ref{sec:comparison},~\ref{sec:morphology_1} 
and~\ref{sec:morphology_2}, all partial derivatives 
with respect to \mbox{time $t$} have been dropped out
in system~\eqref{eq:lumped_model} in such a way that
Eq.~\eqref{Polaron_red} reduces to an algebraic constraint.

The numerical strategy adopted in the present paper is basically
composed of three steps: 
\begin{enumerate}
  \item Linearization
  \item Spatial discretization
  \item Solution of the linear algebraic system
\end{enumerate}

{\it Step (1)}\\
For model linearization,  we adopt a quasi-Newton approach similar to
that used in~\cite{defalco}, where, in the computation of Jacobian matrix 
entries, the dependence of mobilities and polaron pair dissociation rate on the solution
is neglected. 

{\it Step (2)}\\
Similarly to~\cite{defalco}, for the spatial discretization of the sequence of linear systems of PDEs stemming from  {\it Step (1)} we adopt the Galerkin-Finite Element Method (G-FEM) stabilized by means on an Exponential Fitting technique~\cite{bank1998,gatti1998,xu1999,lazarov} 
in order to deal with possibly dominating drift terms in the continuity equations.
A peculiarity of the heterojunction 
model~\eqref{eq:lumped_model} as compared to the \emph{homogenized} model of~\cite{defalco} is the presence of non-trivial interface conditions at the donor-acceptor interface, which are taken care of by means of the \emph{substructuring} techniques described, {\it e.g.}, 
in~\cite{domain,hughes2000} which turn out to be of straightforward implementation in the adopted G-FEM method.

{\it Step (3)}\\
%In the 1- and 2D cases considered thus far, 
To solve the linear algebraic systems arising from problem 
discretization, we employ the  Unsymmetric Multi Frontal method 
implemented in the UMFPACK library~\cite{timdavis}
as on current hardware architectures memory constraints 
are not the main limiting factor and the use of 
a direct sparse solver has the advantage of being
more robust than iterative approaches 
with respect to coefficient matrix conditioning
% \cdf{As on current hardware architectures memory constraints are not usually the main limiting factor,
% we adopt a direct method, namely the  Unsymmetric Multi Frontal method implemented in the UMFPACK libarary~\cite{timdavis},
% to solve linear algebraic systems. Using a direct sparse solver has the advantage of being
% more robust with respect to the coefficient matrix conditioning.}
% because memory consumption is not a major constraint. 
% In future application of the model to the study of 3D morphologies,
% a more computationally efficient approach might be 
% to resort to subdomain iterations commonly employed in
% Domain Decomposition Methods~\cite{domain}.

%%% Local Variables: 
%%% mode: latex
%%% TeX-master: "paper_osc_multiscale"
%%% End: 

\section{Simulation Results}
\label{sec:num_res}
In this section we carry out an extensive computational 
study of the micro and macroscale models introduced in
Sects.~\ref{sec:microscale_model} and~\ref{sec:macroscale_model}.
In Sect.~\ref{sec:comparison_macro_micro}, one-dimensional
transient simulations under different working
conditions are carried out to verify the accuracy 
of the macroscale model with respect
to the microscale system. 
In Sects.~\ref{sec:comparison},~\ref{sec:morphology_1} 
and~\ref{sec:morphology_2}, computations
are performed in steady-state conditions, in order to validate
the macroscale model by comparison with available 
results in the literature.
The numerical schemes of Sect.~\ref{sec:numer_approx} have been 
implemented in \texttt{Octave} using the Octave-Forge package
%\texttt{msh} \cite{msh} and
\texttt{bim} \cite{bim} for matrix assembly.
% Since in this work we focus on steady-state operating conditions, the resulting nonlinear stationary problems are solved with \texttt{fsolve}, nevertheless the code is also able to solve evolutionary problems, using the library \texttt{daspk} \cite{daspk,daspk2}.\\
%The code is available at INDIRIZZO RILASCIO.\\

\subsection{Numerical Validation of the Accuracy of the 
Macroscale Model}
\label{sec:comparison_macro_micro}
In Sects.~\ref{sec:macroscale_model} and~\ref{sec:microscale_model}, 
we have illustrated two models 
of the operating principles of bilayer OSCs 
at two increasing levels of detail, corresponding to 
the macro and micro scales, respectively. 
The two modeling descriptions are expected
to provide a correspondingly more refined level of quality
in the representation of the principal physical phenomena that
govern the functioning of an OSC, at the price, however, of
a substantial increase in implementation complexity and 
computational effort, especially in the case of 
multi-dimensional simulations 
(mesh generation of complex interface morphologies, 
solution of large algebraic systems with possibly badly-balanced
matrices). The natural question that arises at this point of the
discussion is whether the macroscale formulation of 
Sect.~\ref{sec:complete_microscale_model} is capable of 
returning an output picture of the performance of a bilayer OSC
with sufficient accuracy compared to that of the microscale formulation of Sect.~\ref{sec:microscale_model}. By construction of the 
two models, the extent by which the word 
\lq\lq accuracy\rq\rq{} is mathematically identified
cannot certainly refer to a pointwise comparison between 
micro and macroscale solutions (they will certainly look
different!), rather, it should be concerned with 
{\em average quantities} that best represent the overall 
performance of the device. In this respect, the verification test
we are going to carry out later on, is a check of the 
total current per unit area $j_{tot}(t)$
predicted by the micro and macro formulations, where
$j_{tot}(t) = |(\vect{j}(t) + 
\partial{(\varepsilon \vect{E})}/\partial t) \cdot 
\boldsymbol{\nu}|_{\Gamma_{cont}}$, $\Gamma_{cont}$ being 
either of the two contacts $\Gamma_C$ or $\Gamma_A$.
The choice of $j_{tot}$ for model validation 
is due to the fact that the total output current density 
is an easily accessible quantity in experiments, and thus
represents the most significant parameter
for assessing the photoconversion performance of a solar 
cell~\cite{mcgehee}.

\begin{figure}[tb!]
	\centering
	\subfigure[Relative difference in $j_{tot}$.]{
		\includegraphics[width=.475\textwidth]{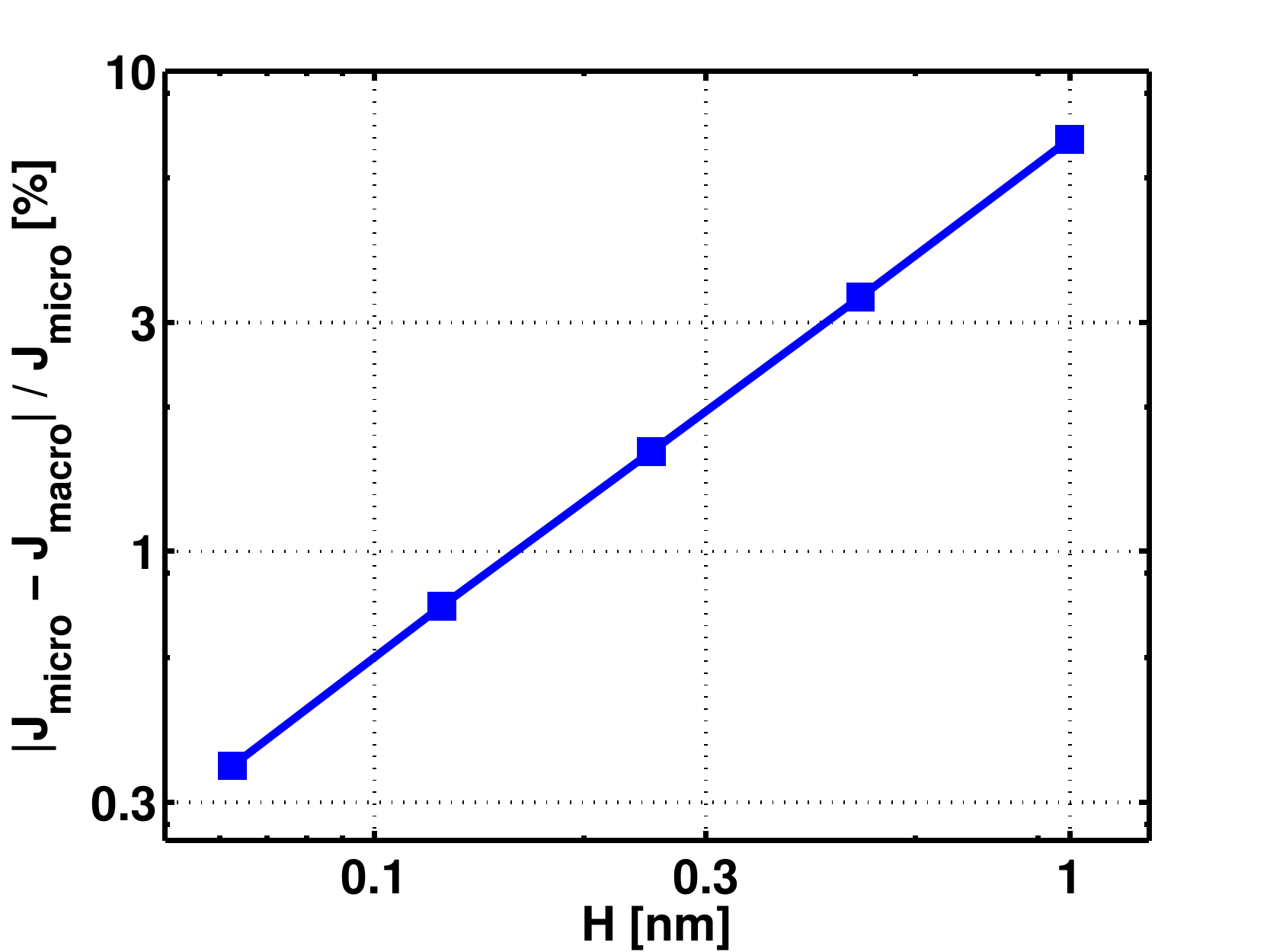}
		\label{fig:comparison_ratio}
	}
	\subfigure[Time evolution of $j_{tot}$.]{
		\includegraphics[width=.475\textwidth]{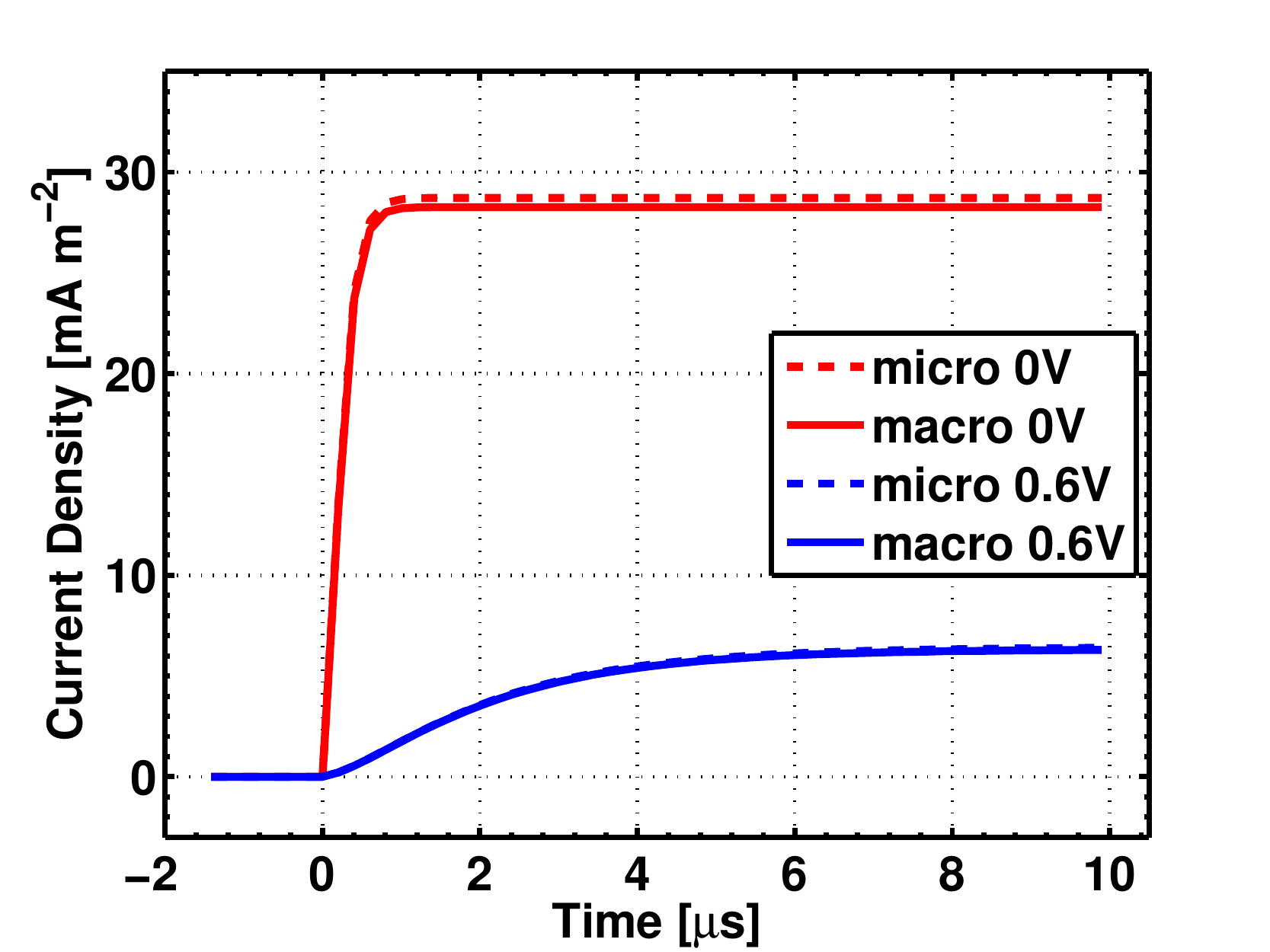}
		\label{fig:comparison_current}
	}
	\caption{Comparison between microscale and macroscale models.} 
	\label{fig:comparison}
\end{figure}

For sake of computational simplicity, we consider a biplanar OSC, 
so that the resulting spatial geometrical description 
can be reduced to a one--dimensional model. 
The total length $L_{cell}$ of the device is equal to $100$\,nm,
with the two regions $\Omega_n$ and $\Omega_p$ occupying
each one half of the cell. All model coefficients are assumed
to be constant quantities, and their values are listed 
in Tab.~\ref{tab:model_param_1d}. 
\begin{table}
\begin{center}
\begin{tabular}{lll}\hline
\textbf{Parameter}							& \textbf{Symbol} 		& \textbf{Numerical value}\\ \hline
Acceptor relative dielectric constant		& $\varepsilon_{r,a}$	& 2.5\\
Donor relative dielectric constant 			& $\varepsilon_{r,d}$	& 2.5\\
Built-in voltage 							& $V_{bi}$				& -0.6$\,$V\\
Temperature 								& $T$					& 298$\,$K\\
Electron mobility & $\mu_{n}$		& $4\cdot 10^{-8}\,\meter^2\,\volt^{-1}\,\second^{-1}$\\
Hole mobility & $\mu_{p}$ 			& $2\cdot 10^{-8}\,\meter^2\,\volt^{-1}\,\second^{-1}$\\
Exciton diffusion coefficient							& $D_e$					& $1\cdot 10^{-7}\,\meter^2\,\second^{-1}$\\
Exciton lifetime							& $\tau_{e}$		& $1\cdot 10^{-9}\,\second$\\
Exciton dissociation time					& $\tau_{diss}$ 		& $1\cdot 10^{-12}\,\second$\\
Polaron pair recombination rate	constant & $k_{rec}$				& $1\cdot 10^{6}\,\second^{-1}$\\
Singlet exciton recombination fraction		& $\eta$				& $0.25$ \\
Polaron pair dissociation rate constant	($V_{appl}=0$\,V) & $k_{diss}$ 			& $1\cdot 10^{7}\,\second^{-1}$\\
Polaron pair dissociation rate constant	($V_{appl}=-V_{bi}=+0.6$\,V)
& $k_{diss}$ 			& $2\cdot 10^{5}\,\second^{-1}$\\
Bimolecular recombination rate constant & $\gamma$ & $1\cdot 10^{-19}\,\meter^3\,\second^{-1}$\\
\multicolumn{3}{l}{Boundary condition parameters for electrons}\\
		& $\kappa_n$&  $0$\\
		& $\alpha_n$ & $1\,\meter\,\second^{-1}$\\
		& $\beta_n$ & $0\,\meter^{-2}\,\second^{-1}$\\[2mm]
\multicolumn{3}{l}{Boundary condition parameters for holes}\\
		& $\kappa	_p$&  $0$\\
		& $\alpha_p$ & $1\,\meter\,\second^{-1}$\\
		& $\beta_p$ & $0\,\meter^{-2}\,\second^{-1}$\\[2mm]
\hline
\end{tabular}

\caption{Model parameter values used in the simulations
of Sect.~\ref{sec:comparison_macro_micro}.}
\label{tab:model_param_1d}
\end{center}
\end{table}
We start by simulating the cell turn-on transient 
at short circuit condition ($V_{appl} = 0\,\volt$), 
which corresponds to computing the device response 
to an abrupt variation at $t=0$ in the photon absorption rate 
$Q$ from zero to $10^{25} \, \meter^{-3}\,\second^{-1}$.
The simulation time interval 
is taken wide enough for the device to reach stationary conditions.
In Fig.~\ref{fig:comparison_ratio} the relative discrepancy between
the stationary values of  $j_{tot}$
computed with the two methods is reported for several values of the
interface width parameter $H$. Results allow us 
to conclude that in the physically relevant range of 
variation of $H$, the relative discrepancy between the 
micro and macroscale models remains consistently below 10\%
and that, as expected, the predictions of the two models tend 
to become undistinguishable as $H$ tends to zero.
In Fig.~\ref{fig:comparison_current} we set $H = 0.25\,\nano\meter$
and we show the time evolution 
of the total current per unit area in the two biasing conditions 
$V_{appl} = 0\,\volt$ and $V_{appl} = -V_{bi} = 0.6\,\volt$ this latter 
being the ``flat band'' voltage.
Accordingly to Fig.~\ref{fig:comparison_ratio}, in both regimes
the curves almost coincide over the whole simulation time interval.

\subsection{Model Validation through Comparison with Existing 
Simulation Data}\label{sec:comparison}
In this section, we aim to compare the predictions of our macroscale model to those 
of \cite{williams_th,williams2008} 
and to investigate the impact of the model 
for $k_{diss}$ proposed in 
Sect.~\ref{sec:k_diss_model}
on the simulated device performance.
We consider the same device as in~ \cite{williams_th,williams2008,CogliatiPorroThesis2010}  
where the  acceptor and donor materials are F8BT and PFB, 
respectively.
The values of the model parameters are listed in Table~\ref{tab:model_param}.

% ,
% except for the exciton transition characteristic time $\tau_{diss}$ 
% that has been set equal to $10^{-12}\,\second$ according to~\cite{?????}.
% since that phenomenon is known to be very fast if compared to the other occurring in the cell.

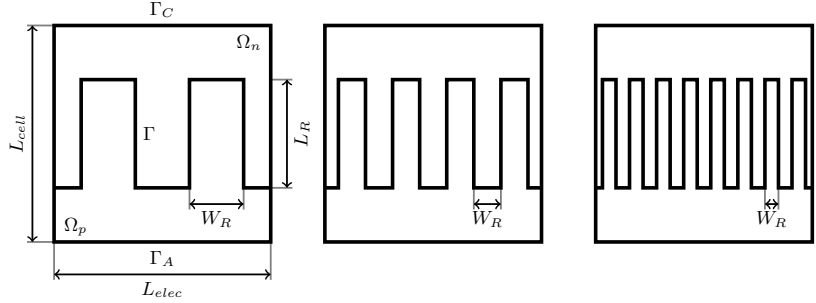
\begin{figure}[!htb]
	\centering	
	\scalebox{.8}{\begin{tikzpicture}[scale=1.8]
%  	\draw [<->,thick, line width=0.4mm] (0,2.3) node (yaxis) [left] {$y$}
%       		 |- (2.3,0) node (xaxis) [below] {$x$};

	\draw[line width=0.3mm, gray] (0, -0.35) -- (0, 0);
	\draw[line width=0.3mm, gray] (2, -0.35) -- (2, 0);
	\draw[line width=0.3mm, gray] (-0.25, 0) -- (0, 0);
	\draw[line width=0.3mm, gray] (-0.25, 2) -- (0, 2);
	\draw[line width=0.3mm, gray] (2, 0.5) -- (2.2, 0.5);
	\draw[line width=0.3mm, gray] (2, 1.5) -- (2.2, 1.5);	
	\draw[line width=0.3mm, gray] (1.25, 0.3) -- (1.25, 0.5);
	\draw[line width=0.3mm, gray] (1.75, 0.3) -- (1.75, 0.5);
	
	\draw[line width=0.6mm] (0,0) rectangle (2,2);
	\draw[line width=0.6mm] (0, 0.5) 	-- (0.25, 0.5)
									 	-- (0.25, 1.5)
									 	-- (0.75, 1.5)
										-- (0.75, 0.5)
										-- (1.25, 0.5)
										-- (1.25, 1.5)
										-- (1.75, 1.5)
										-- (1.75, 0.5)
										-- (2, 0.5);
										
	\draw[<->,thick, line width=0.3mm] (0, -0.3) -- (2, -0.3) node[sloped,midway,below] {$L_{elec}$};
	\draw[<->,thick, line width=0.3mm] (-0.2, 0) -- (-0.2, 2) node[sloped,midway,above] {$L_{cell}$};
	\draw[<->,thick, line width=0.3mm] (2.15, 0.5) -- (2.15, 1.5) node[sloped,midway,below] {$L_R$};
	\draw[<->,thick, line width=0.3mm] (1.25, 0.35) -- (1.75, 0.35) node[sloped,midway,below] {$W_R$};

    \draw (1.8, 1.84) node {$\Omega_{n}$};
	\draw (0.2, 0.14) node {$\Omega_{p}$};
	\draw (0.88, 1) node {$\Gamma$};
	\draw (1, -0.15) node {$\Gamma_{A}$};
	\draw (1, 2.15) node {$\Gamma_{C}$};
%	\draw (-0.16, 1) node {$\Gamma_{N}$};
%	\draw (2.17, 1) node {$\Gamma_{N}$};

% Second geometry
	\draw[line width=0.3mm, gray] (2.5+1.375, 0.3) -- (2.5+1.375, 0.5);
	\draw[line width=0.3mm, gray] (2.5+1.625, 0.3) -- (2.5+1.625, 0.5);

	\draw[line width=0.6mm] (2.5+0,0) rectangle (2.5+2,2);
	
	\draw[line width=0.6mm] (2.5+0, 0.5)	-- (2.5+0.125, 0.5)
											-- (2.5+0.125, 1.5)
											-- (2.5+0.375, 1.5)
											-- (2.5+0.375, 0.5)
											-- (2.5+0.625, 0.5)
											-- (2.5+0.625, 1.5)
											-- (2.5+0.875, 1.5)
											-- (2.5+0.875, 0.5)
											-- (2.5+1.125, 0.5)
											-- (2.5+1.125, 1.5)
											-- (2.5+1.375, 1.5)
											-- (2.5+1.375, 0.5)
											-- (2.5+1.625, 0.5)
											-- (2.5+1.625, 1.5)
											-- (2.5+1.875, 1.5)
											-- (2.5+1.875, 0.5)
											-- (2.5+2, 0.5);

	\draw[<->,thick, line width=0.3mm] (2.5+1.375, 0.35) -- (2.5+1.625, 0.35) node[sloped,midway,below] {$W_R$};	

% Third geometry
	\draw[line width=0.3mm, gray] (6+1.125/2, 0.3) -- (6+1.125/2, 0.5);
	\draw[line width=0.3mm, gray] (6+1.375/2, 0.3) -- (6+1.375/2, 0.5);	
	
	\draw[line width=0.6mm] (5+0,0) rectangle (5+2,2);
	
	\draw[line width=0.6mm] (5+0, 0.5)	-- (5+0.125/2, 0.5)
										-- (5+0.125/2, 1.5)
										-- (5+0.375/2, 1.5)
										-- (5+0.375/2, 0.5)
										-- (5+0.625/2, 0.5)
										-- (5+0.625/2, 1.5)
										-- (5+0.875/2, 1.5)
										-- (5+0.875/2, 0.5)
										-- (5+1.125/2, 0.5)
										-- (5+1.125/2, 1.5)
										-- (5+1.375/2, 1.5)
										-- (5+1.375/2, 0.5)
										-- (5+1.625/2, 0.5)
										-- (5+1.625/2, 1.5)
										-- (5+1.875/2, 1.5)
										-- (5+1.875/2, 0.5)
										-- (5+2.125/2, 0.5)

										-- (6+0.125/2, 1.5)
										-- (6+0.375/2, 1.5)
										-- (6+0.375/2, 0.5)
										-- (6+0.625/2, 0.5)
										-- (6+0.625/2, 1.5)
										-- (6+0.875/2, 1.5)
										-- (6+0.875/2, 0.5)
										-- (6+1.125/2, 0.5)
										-- (6+1.125/2, 1.5)
										-- (6+1.375/2, 1.5)
										-- (6+1.375/2, 0.5)
										-- (6+1.625/2, 0.5)
										-- (6+1.625/2, 1.5)
										-- (6+1.875/2, 1.5)
										-- (6+1.875/2, 0.5)
										-- (6+2/2, 0.5);
									
	\draw[<->,thick, line width=0.3mm] (6+1.125/2, 0.35) -- (6+1.375/2, 0.35) node[sloped,midway,below] {$W_R$};	

\end{tikzpicture}}
	\vspace{-0.3cm}	
	\caption{Internal morphology with rod-shaped donor-acceptor interface.}
	\label{fig:finger_geo}
%	\vspace{-1cm}
\end{figure}

The device morphology, shown in Fig.~\ref{fig:finger_geo}, 
is an interpenetrating rod-shaped structure of donor and acceptor materials
with $L_{cell} = 150\,\nano\meter$, $L_{elec} = 50\,\nano\meter$, $L_R = 79\,\nano\meter$ and  $W_R = 6.25\,\nano\meter$.
Throughout this section, we denote by $y$ the direction between the two electrodes $\Gamma_C$ and $\Gamma_A$.

In \cite{williams_th,williams2008} an optical model has been used to
determine the exciton generation term $Q$. Here, instead, we follow a
simpler approach by considering $Q$ to be constant in the entire
device structure  
and equal to the value obtained averaging the result
in~\cite{williams_th,williams2008}. 
The choice $\kappa_n = \kappa_p = 0$ corresponds to enforcing Dirichlet
boundary conditions for the carrier densities at the contacts and
amounts to neglecting the dependence of charge
injection on the electric field and assuming an infinite recombination rate at the contacts.
% Dirichlet boundary 
% conditions for $n$ and $p$ are enforced at the electrodes, unlike 
% 
% are simplified with respect to those mentioned
% in Sect.~\ref{sec:reduced_model}, Dirichlet type boundary
% conditions for charge carriers are enforced at the electrodes too,
% obtained by assuming that charge injection is negligible there since
% electrons and holes experience high energy barriers. 
%Following \cite{williams_th,williams2008} we assume the device to be
%irradiated with monochromatic light with wavelength in the range
%459-460$\,\nano\meter$ and define $1\,$Sun as the power density of
%that range in the AM 1.5 spectrum.

\begin{figure}[!tbh]
	\centering
	\subfigure[$Q = 1.53 \cdot 10^{23}\,\meter^{-3}\,\second^{-1}$]{
	\includegraphics[width=0.475\textwidth]{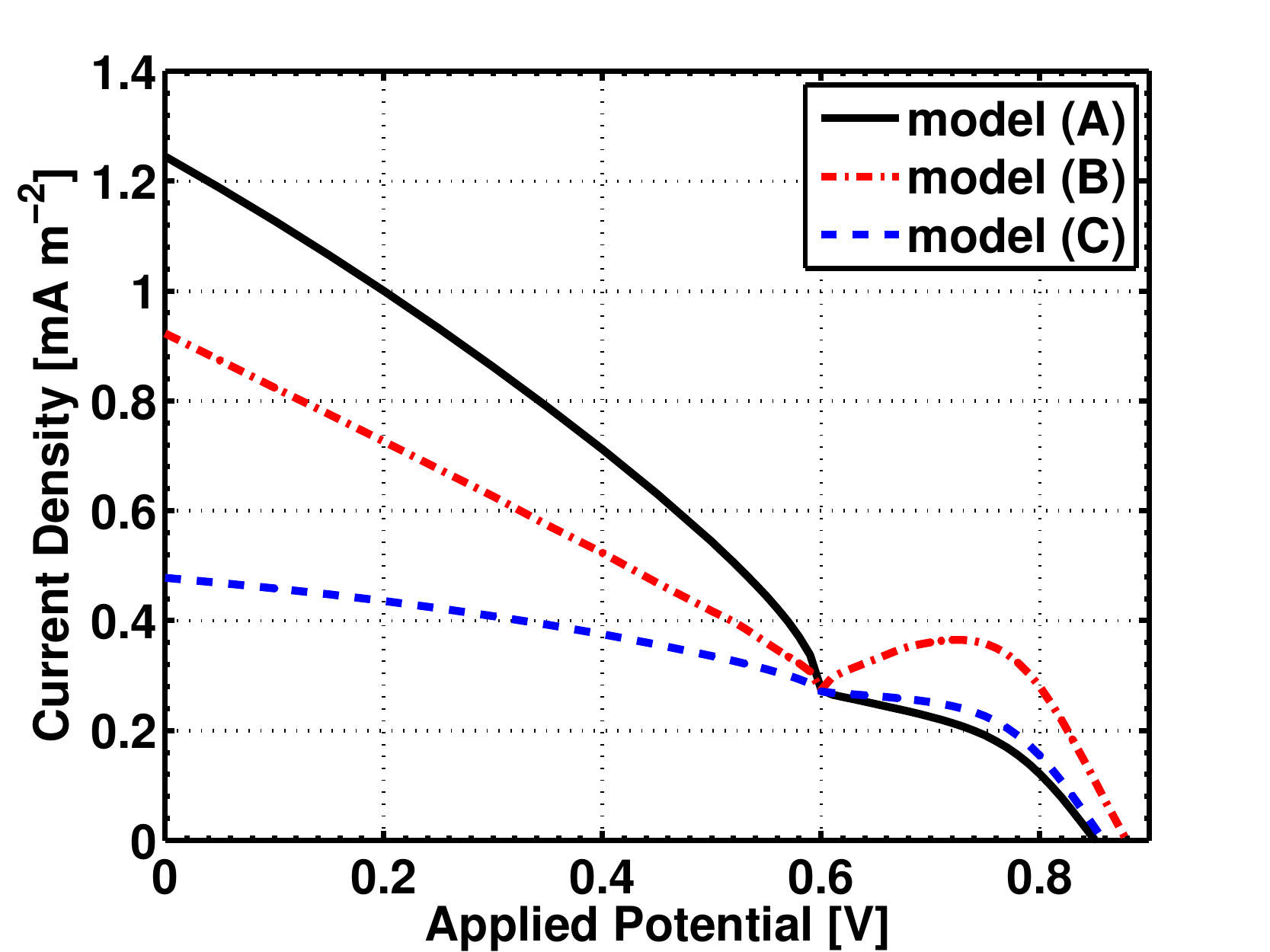}
	\label{fig:iv_comparison_left}
	}
	\subfigure[$Q = 1.53 \cdot 10^{25}\,\meter^{-3}\,\second^{-1}$]{
	\includegraphics[width=0.475\textwidth]{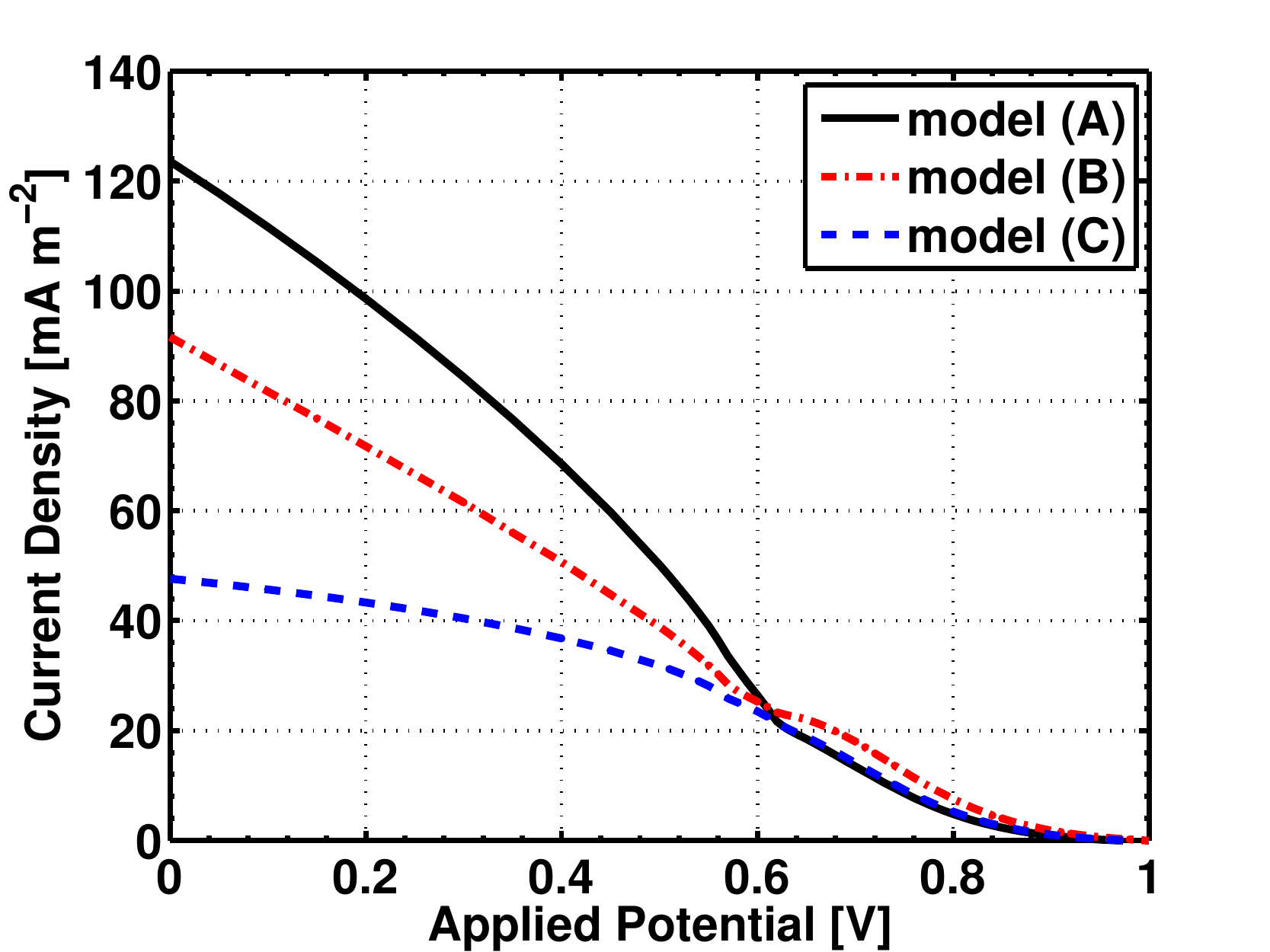}
	\label{fig:iv_comparison_right}
	}
	\caption{Comparison of the current-voltage characteristic lines 
	with two different
	values for the exciton generation rate. Solid line curve refers to 
	model proposed in~\cite{williams_th}, dash-dotted line refers
        to model~\eqref{eq:kdissmodel2} while dashed line
	refers to model~\eqref{eq:kdissmodel3}.}
	\label{fig:iv_comparison}
\end{figure}

Fig.~\ref{fig:iv_comparison_left}
shows the current density-voltage characteristics in the 
case of an exciton generation rate $Q = 1.53 \cdot 10^{23}\,\meter^{-3}\,\second^{-1}$.
\begin{table}
\begin{center}
\begin{tabular}{lll}\hline
\textbf{Parameter}							& \textbf{Symbol} 		& \textbf{Numerical value}\\ \hline
Acceptor relative dielectric constant		& $\varepsilon_{r,a}$	& 4\\
Donor relative dielectric constant 			& $\varepsilon_{r,d}$	& 4\\
%Built-in voltage 							& $V_{bi}$				& 0.6$\,$V\\
Temperature 								& $T$					& 298$\,$K\\[2mm]
\multicolumn{3}{l}{Poole-Frenkel mobility model parameters for electrons~\cite{barker2003}}\\
& $\mu_{n}(0)$		& $3\cdot 10^{-10}\,\meter^2\,\volt^{-1}\,\second^{-1}$\\
& $\gamma_a$			& $1.55\cdot 10^{-3}\,\volt^{-1/2}\,\meter^{1/2}$\\[2mm]
\multicolumn{3}{l}{Poole-Frenkel mobility model parameters for holes~\cite{barker2003}}\\
& $\mu_{p}(0)$ 			& $1\cdot 10^{-10}\,\meter^2\,\volt^{-1}\,\second^{-1}$\\
& $\gamma_d$			& $3\cdot 10^{-4}\,\volt^{-1/2}\,\meter^{1/2}$\\[2mm]
%Cathode injection barrier height			& $\phi_{bn}^c$			& 0.5$\,$eV\\
%Anode injection barrier height				& $\phi_{bp}^a$			& 0.5$\,$eV\\
Exciton diffusion coefficient							& $D_e$					& $1\cdot 10^{-7}\,\meter^2\,\second^{-1}$\\
Exciton lifetime							& $\tau_{e}$		& $1\cdot 10^{-9}\,\second$\\
Exciton dissociation time					& $\tau_{diss}$ 		& $1\cdot 10^{-12}\,\second$\\
Polaron pair recombination rate	constant			& $k_{rec}$				& $1\cdot 10^{6}\,\second^{-1}$\\
Singlet exciton recombination fraction		& $\eta$				& $0.25$ \\
Polaron pair zero-field dissociation rate constant	& $k_{diss}(0)$ 			& $1\cdot 10^{5}\,\second^{-1}$\\
Interface width  					& $2H$					& $2\cdot 10^{-9}\,\meter$\\[2mm]
%Conduction band states density				& $N_C$					& $10^{27}\,\meter^{-3}$\\
%Valence band states density					& $N_D$					& $10^{27}\,\meter^{-3}$\\
\multicolumn{3}{l}{Boundary condition parameters for electrons~\cite{defalco}}\\
		& $\kappa_n$&  $0$\\
		& $\alpha_n$ & $1\,\meter\,\second^{-1}$\\
		& $\beta_n$ & $3.4995\cdot 10^{18}\,\meter^{-2}\,\second^{-1}$\\[2mm]
\multicolumn{3}{l}{Boundary condition parameters for holes~\cite{defalco}}\\
		& $\kappa	_p$&  $0$\\
		& $\alpha_p$ & $1\,\meter\,\second^{-1}$\\
		& $\beta_p$ & $3.4995\cdot 10^{18}\,\meter^{-2}\,\second^{-1}$\\[2mm]
%Space scaling								& $\overline{L}$		& $150\cdot 10^{-9}\,\meter$\\
%Carrier scaling								& $\overline{n}$		& $1\cdot 10^{22}\,\meter^{-3}$\\
%Exciton scaling								& $\overline{X}$		& $1\cdot 10^{12}\,\meter^{-2}$\\
\hline
\end{tabular}

\caption{Model parameter values used in the simulations of Sects.~\ref{sec:comparison}, 
\ref{sec:morphology_1} and~\ref{sec:morphology_2}.}
\label{tab:model_param}
\end{center}
\end{table}
The three curves correspond to the use of three different 
expressions for the polaron dissociation rate constant $k_{diss}$,
identified as follows: (A) the model proposed 
in~\cite{williams_th,williams2008} with
\mbox{$\overline{E}_y = {\left|\Gamma\right|}^{-1}\!\!\int_{\Gamma}E_y\,dx$}
as the driving parameter for polaron pair dissociation (solid line); (B) 
the model~\eqref{eq:kdissmodel2} (dash-dotted line); (C)
the model~\eqref{eq:kdissmodel3} (dashed line).
The result computed using model (A)
is in excellent agreement with that of Fig.~7(right) in~\cite{williams2008} despite the above mentioned modeling differences.
Model (A)
does not account for the orientation of the electric field with respect to the donor-acceptor interface and is expected to 
overestimate dissociation in the case where 
$\mathbf{E} \cdot \boldsymbol{\nu}_{\Gamma} \simeq 0$. This is confirmed by
the curve for model (B). As a matter of fact, in this
case all the dissociation directions are assumed to be equally likely and the 
%when model \eqref{eq:kdissmodel2} is used, accounting for both normal
%and tangential components of the electric field.
computed output current density before flat-band condition
occurs ($V_{appl} \leq 0.6\,\volt$) is smaller than predicted 
by the solid line curve. 
For $V_{appl} > 0.6\,\volt$ 
the computed current-voltage characteristic exhibits a non-monotonic
behavior. This latter behavior is not observed in any of 
the experimental measurements we are aware of, 
and is most probably  
to be ascribed to a too important contribution of the 
tangential component of the electric field $E_t$ that
leads~\eqref{eq:kdissmodel2} to overestimate polaron 
dissociation at the material interface.
If, instead, model (C) is used, 
the obtained output current density characteristics is 
the dashed line in Fig.~\ref{fig:iv_comparison_left}.
We observe a smoother trend than in previous cases for all applied
voltages,  
and close to short circuit condition we note that the current density
is further reduced since dissociation is assumed to occur only in the
normal direction and on a significant portion of the interface $E_n$
is almost vanishing. 
In all the considered cases, the nonsmooth behavior 
at flat band condition ($V_{appl} = 0.6\,\volt$) is to be ascribed to
the discontinuity of $\partial \beta/\partial z$ at $z=0$ 
in~\eqref{eq:enhancement}. 

Fig.~\ref{fig:iv_comparison_right} shows the results of the same
analysis 
as above in the case of an exciton generation rate $Q = 1.53 \cdot
10^{25}\,\meter^{-3}\,\second^{-1}$. The shape of the characteristics
is very similar to those with low light up to a scaling factor of
about 100, this suggesting a linearity between the output current
density and the illumination intensity. Notice the absence of the bump
for $V_{appl} > 0.6\,\volt$ in the 
case of model (B). This is a consequence of the increased magnitude of
the charge carrier densities compared to the previously considered
illumination that in turn determines stronger Coulomb attraction
forces and hence more recombination phenomena. With reduced
attractions, instead, charge carriers have more chances to escape from
the interface following concentration gradients. 

\begin{figure}[!t]
	\centering
	\includegraphics[width=0.44\textwidth]{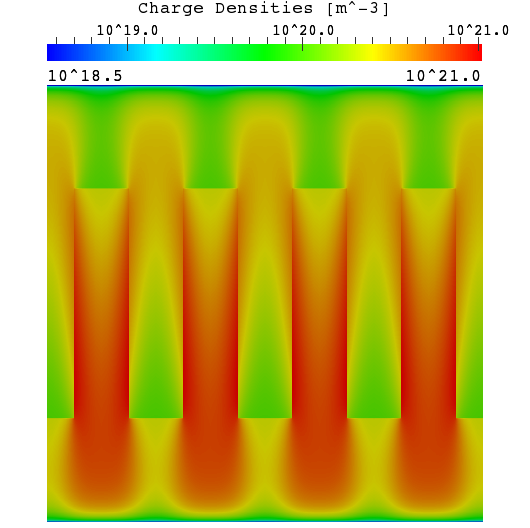}
	\includegraphics[width=0.44\textwidth]{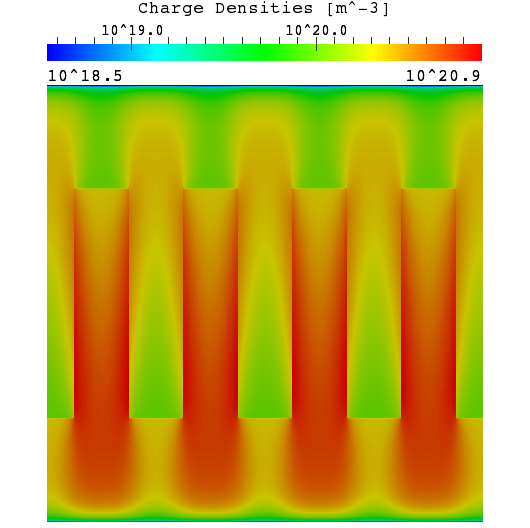} \\
	\includegraphics[width=0.44\textwidth]{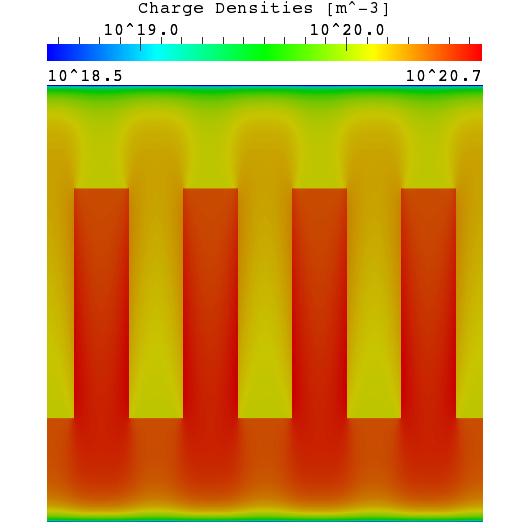}
	\caption{Charge carrier densities $[\meter^{-3}]$ at short
          circuit condition with 
	\mbox{$Q = 1.53 \cdot 10^{25}\,\meter^{-3}\,\second^{-1}$}
        using models (A) (left), (B) (right) and 
	(C) (bottom), respectively.}
	\label{fig:densities}
\end{figure}

In Fig.~\ref{fig:densities} 
we show the charge carrier densities in a device with geometrical data
set to $L_{cell} = 150\,\nano\meter$, $L_{elec} = 440\,\nano\meter$,
$L_R = 79\,\nano\meter$ and $W_R = 55\,\nano\meter$, at short circuit
condition with exciton generation rate 
\mbox{$Q = 1.53 \cdot 10^{25}\,\meter^{-3}\,\second^{-1}$}. 
We first observe that computed charge carrier distributions in
Fig.~\ref{fig:densities}(left) are in very good agreement with those
of Fig.~3(i) in~\cite{williams_th}
%
% except near the electrodes where 
%different boundary conditions are applied, 
%\begin{center}
%\nuovo{QUI IL REFEREE 3 HA OBIETTATO QUALCOSA SULLE BCS. 
%SISTEMIAMO LA QUESTIONE NEL MODO PIU` SEMPLICE}
%\end{center}
%
and show the same peaks
close to the vertical sides of the donor-acceptor interface. It is
interesting to notice  
that the total number of holes in the donor material is higher than
the number of electrons in the acceptor material because of the
significantly different values of their respective
mobilities. Negative charges can move through the device faster to be
finally extracted at the cathode so that an  
overall positive charge builds-up in the device.
The charge densities computed using models (B) and (C) exhibit a
qualitatively similar profile with a gradual reduction of the
magnitude compared to 
the result of model (A). This behavior is completely consistent 
with the previous analysis of the current-voltage characteristics 
predicted by the three models of $k_{diss}$. 

\begin{figure}[!t]
\centering
\subfigure[Open circuit
voltage]{\includegraphics[width=0.475\textwidth]{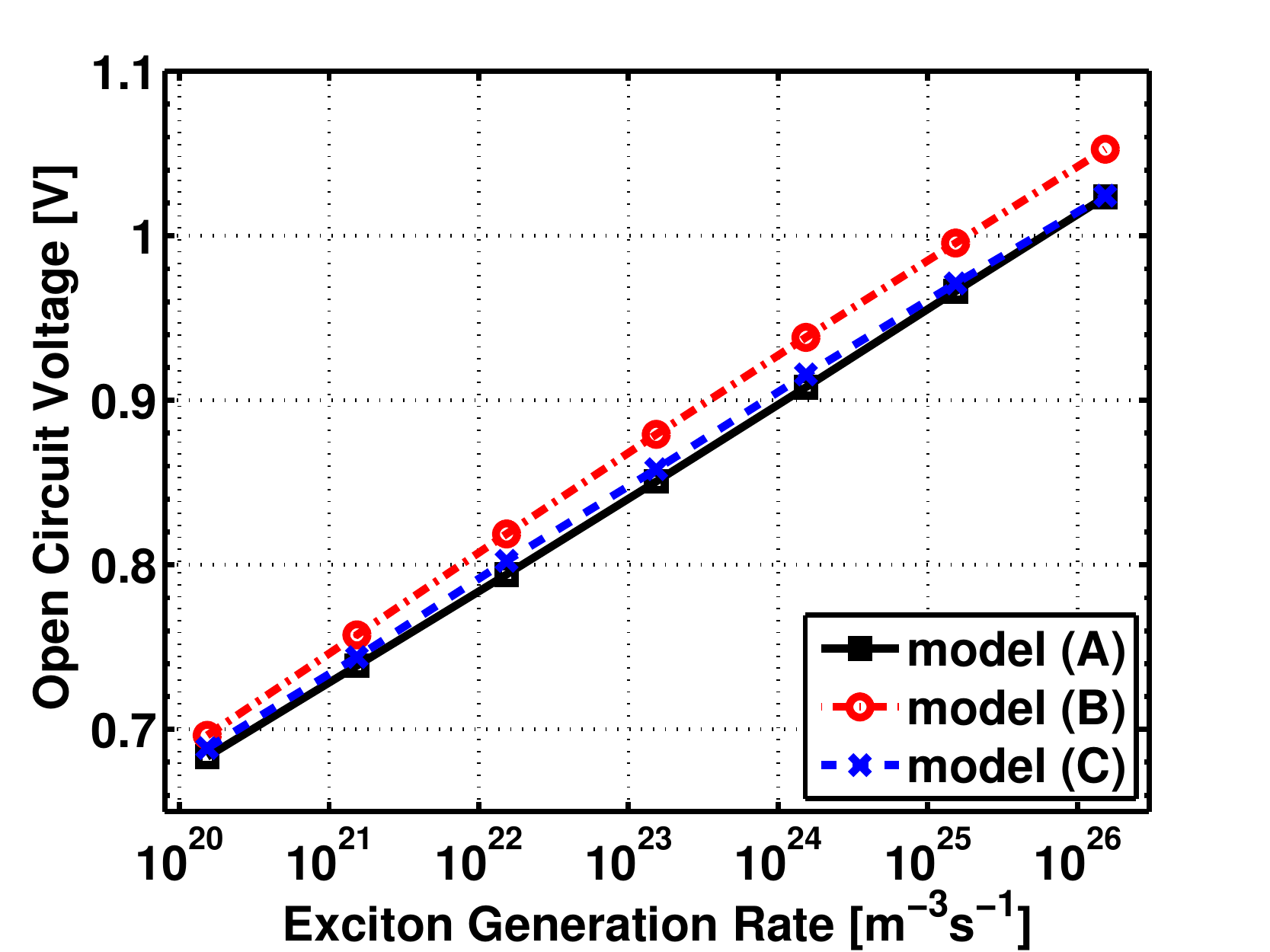} \label{fig:voc_isc_voc}} 
\subfigure[Short circuit current
density]{\includegraphics[width=0.475\textwidth]{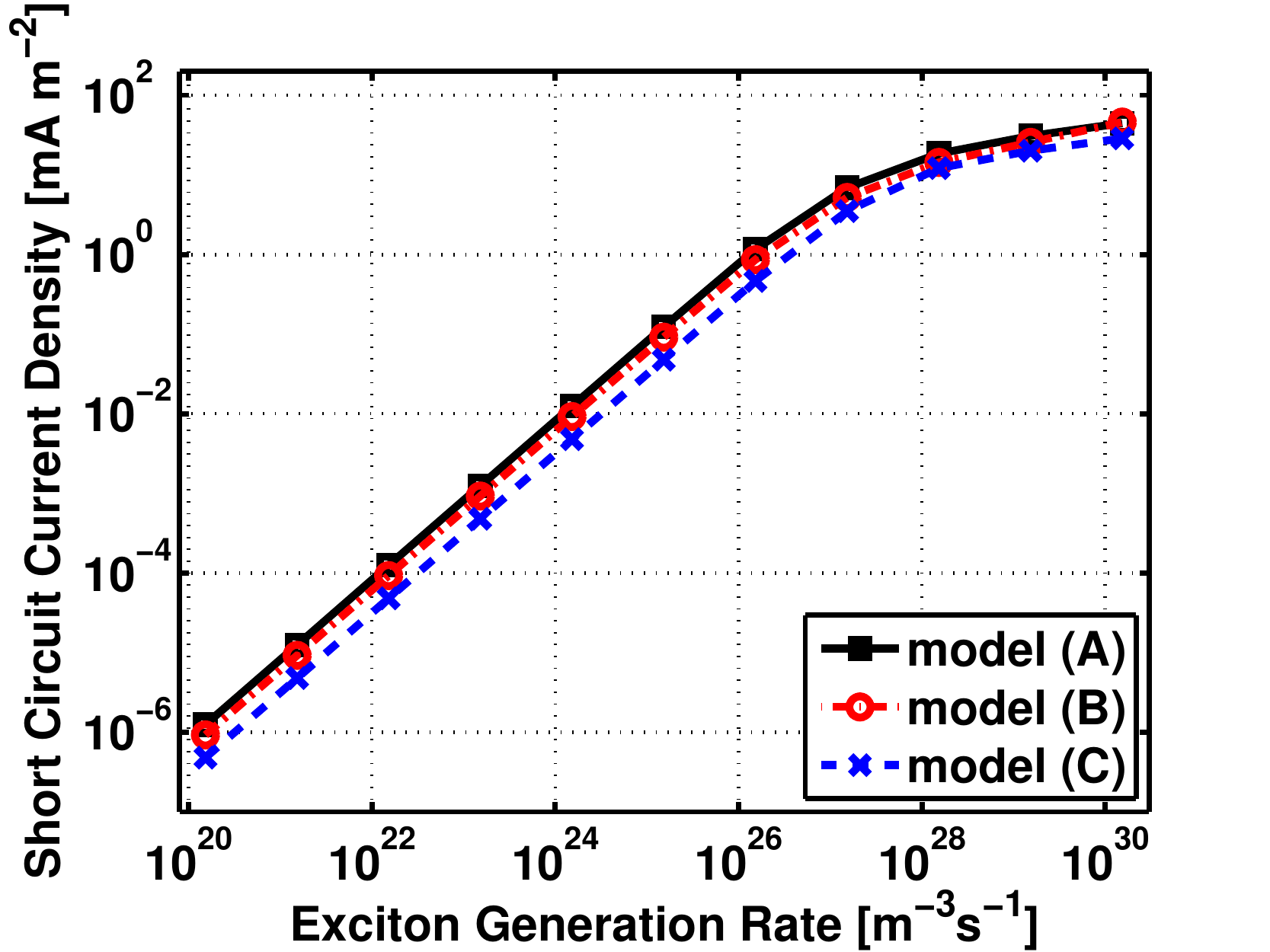}\label{fig:voc_isc_isc}} \caption
{Open circuit voltage and short circuit current density as functions of the exciton generation rate.}\label{fig:voc_isc}
\end{figure}

We conclude this preliminary validation of model \eqref{eq:lumped_model}
by illustrating in Fig.~\ref{fig:voc_isc} the open circuit voltage $V_{oc}$ and short circuit current density $J_{sc}$ of a device with the same characteristics as in the previous set of simulations for values of exciton generation rate in the range from $1.53 \cdot 10^{20}$ to $1.53 \cdot 10^{30}\,\meter^{-3}\,\second^{-1}$. 
Fig.~\ref{fig:voc_isc_voc} is in excellent agreement 
with Fig.~6(right) of~\cite{williams2008}, and
indicates that models (A), (B) and (C) predict a linear behavior of $V_{oc}$ with respect to the logarithm of the exciton generation rate, as already pointed out in \cite{williams_th,williams2008}.
Fig.~\ref{fig:voc_isc_isc} illustrates the current density $J_{sc}$ that can be extracted from the device at short circuit condition. The log-scale plot indicates that $J_{sc}$ increases linearly in a wide range of illumination regimes until values of the order of $10^{28}\,\meter^{-3}\,\second^{-1}$. With more intense irradiation the increase becomes sublinear, suggesting that saturation of the device occurs due to more relevant excitonic and electron-hole recombination phenomena which in turn are a consequence of the increased densities.

\subsection{The Role of Interface Morphology}
\label{sec:morphology_1}
In this section, we aim to investigate the role of 
interface configuration in affecting the OSC performance.
%in order to obtain indications toward an optimal device design.
Referring to Fig.~\ref{fig:finger_geo}, we 
set $L_{cell} = L_{elec} = 150\,\nano\meter$ and 
$L_{R}= 75\,\nano\meter$, and we analyze the importance of interfacial length
by considering devices with an increasing density of interpenetrating structures,
starting from a biplanar device and then taking decreasing values for the rod width $W_R$. Model parameters are the same as in the previous simulations and the exciton generation rate is $Q = 1.53 \cdot 10^{25}\,\meter^{-3}\,\second^{-1}$.

\begin{figure}[hbt!]
	\centering
	\includegraphics[width=0.75\textwidth]{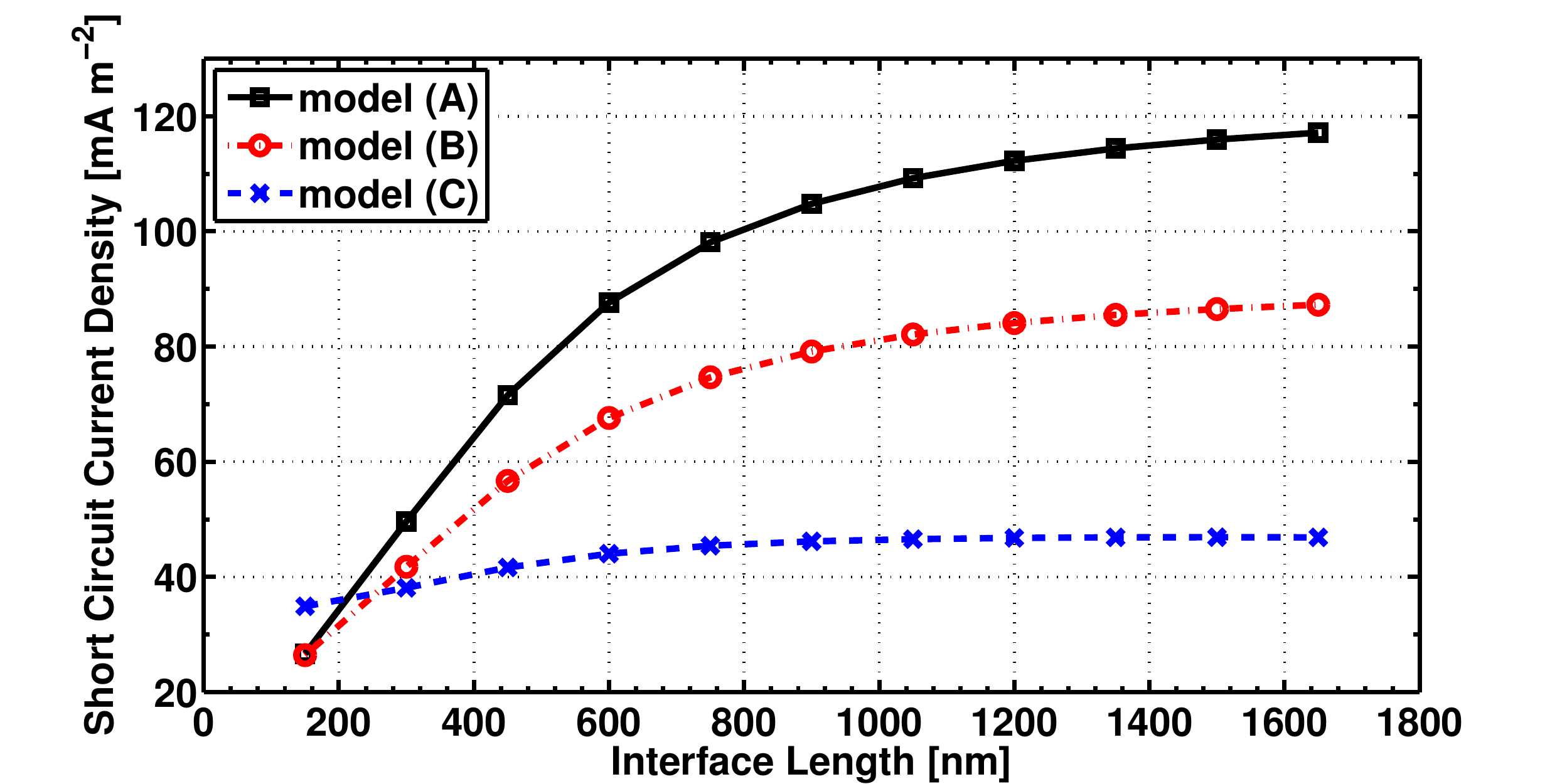}
	\caption{Short circuit current density as a function of interface length.}
	\label{fig:isc_vs_length}
\end{figure}

Fig.~\ref{fig:isc_vs_length} illustrates the computed short circuit current density as a function of the interfacial length for the various polaron pair dissociation rate models we previously considered in this section. In all cases, current saturation is 
predicted for high densities of nanostructures due to the depletion of excitons in the interface area that in turn is a consequence of the abundance of dissociation sites. Computed saturation levels greatly differ among the 
three choices of the model for $k_{diss}$, in accordance with the 
analysis of Sect.~\ref{sec:comparison}.
Fig.~\ref{fig:isc_vs_length} also shows that when a biplanar device is considered, using model (C) a higher short circuit current density is obtained compared to the other approaches. An explanation of this result is 
that the electric field in this case is actually vertically directed 
and this fact, combined with the assumption that dissociation occurs
only in the normal direction, brings to overestimate its rate
(cf. the solid lines in Fig.~\ref{fig:kdissvsE}).
Qualitatively similar results have been obtained in~\cite{buxton2007,williams_th,williams2008}.

%\noindent In the works previously available in literature several words have been spent on analyzing how the length of the interface determines the amount of current that can be extracted from a device with ordered internal nanostructure, see \cite{buxton2007,williams_th,williams2008}. Although the models used therein have different features, the common result has been to show that rising the density of interfacial sites by recurring to more complex morphologies, it is theoretically possible to obtain better charge dissociation even if an increasing saturation phenomenon occurs when the geometrical characteristic dimension comes to the order of nanometer.\\

Also the orientation of the interface is expected to play a role in determining device operation and the following set of simulations aims to investigate this issue. This is a distinctive feature of our model that, to our knowledge, 
has not been treated in previous works. 
For a proper analysis, we allow the orientation of the donor-acceptor interface to change while its overall length remains almost constant, 
in order to single out the effect of the former and analyze it. 

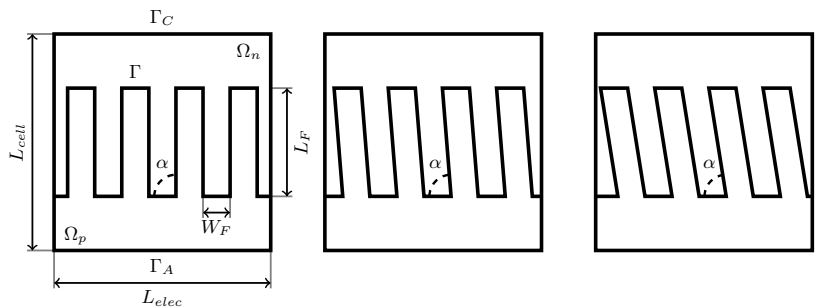
\begin{figure}[h!]
	\centering
	\scalebox{.8}{\begin{tikzpicture}[scale=1.8]
%  	\draw [<->,thick, line width=0.4mm] (0,2.3) node (yaxis) [left] {$y$}
%       		 |- (2.3,0) node (xaxis) [below] {$x$};

	\draw[line width=0.3mm, gray] (0, -0.35) -- (0, 0);
	\draw[line width=0.3mm, gray] (2, -0.35) -- (2, 0);
	\draw[line width=0.3mm, gray] (-0.25, 0) -- (0, 0);
	\draw[line width=0.3mm, gray] (-0.25, 2) -- (0, 2);
	\draw[line width=0.3mm, gray] (2, 0.5) -- (2.2, 0.5);
	\draw[line width=0.3mm, gray] (2, 1.5) -- (2.2, 1.5);	
	\draw[line width=0.3mm, gray] (1.375, 0.3) -- (1.375, 0.5);
	\draw[line width=0.3mm, gray] (1.625, 0.3) -- (1.625, 0.5);
	
	\draw[line width=0.6mm] (0,0) rectangle (2,2);
	
	\draw[line width=0.6mm] (0, 0.5) 	-- (0.125, 0.5)
										-- (0.125, 1.5)
										-- (0.375, 1.5)
										-- (0.375, 0.5)
										-- (0.625, 0.5)
										-- (0.625, 1.5)
										-- (0.875, 1.5)
										-- (0.875, 0.5)
										-- (1.125, 0.5)
										-- (1.125, 1.5)
										-- (1.375, 1.5)
										-- (1.375, 0.5)
										-- (1.625, 0.5)
										-- (1.625, 1.5)
										-- (1.875, 1.5)
										-- (1.875, 0.5)
										-- (2, 0.5);
										
    \draw[dashed, line width=0.4mm] (1.125-0.2, 0.5) arc (180:90:0.2);
%    \draw[dashed, line width=0.6mm] (0.625, 0.5) -- (0.875, 0.5);
	\draw[<->,thick, line width=0.3mm] (0, -0.3) -- (2, -0.3) node[sloped,midway,below] {$L_{elec}$};
	\draw[<->,thick, line width=0.3mm] (-0.2, 0) -- (-0.2, 2) node[sloped,midway,above] {$L_{cell}$};
	\draw[<->,thick, line width=0.3mm] (2.15, 0.5) -- (2.15, 1.5) node[sloped,midway,below] {$L_F$};
	\draw[<->,thick, line width=0.3mm] (1.375, 0.35) -- (1.625, 0.35) node[sloped,midway,below] {$W_F$};
	
    \draw (1.8, 1.84) node {$\Omega_{n}$};
	\draw (0.2, 0.14) node {$\Omega_{p}$};
	\draw (0.75, 1.65) node {$\Gamma$};
	\draw (1, -0.15) node {$\Gamma_{A}$};
	\draw (1, 2.15) node {$\Gamma_{C}$};
%	\draw (-0.16, 1) node {$\Gamma_{N}$};
%	\draw (2.17, 1) node {$\Gamma_{N}$};
	\draw (1, 0.8) node {$\alpha$};

% Second geometry	
	\draw[line width=0.6mm] (2.5+0,0) rectangle (2.5+2,2);
	
	\draw[line width=0.6mm] (2.5+0, 0.5) 	-- (2.5+0.125+0.04, 0.5)
											-- (2.5+0.125-0.04, 1.5)
											-- (2.5+0.375-0.04, 1.5)
											-- (2.5+0.375+0.04, 0.5)
											-- (2.5+0.625+0.04, 0.5)
											-- (2.5+0.625-0.04, 1.5)
											-- (2.5+0.875-0.04, 1.5)
											-- (2.5+0.875+0.04, 0.5)
											-- (2.5+1.125+0.04, 0.5)
											-- (2.5+1.125-0.04, 1.5)
											-- (2.5+1.375-0.04, 1.5)
											-- (2.5+1.375+0.04, 0.5)
											-- (2.5+1.625+0.04, 0.5)
											-- (2.5+1.625-0.04, 1.5)
											-- (2.5+1.875-0.04, 1.5)
											-- (2.5+1.875+0.04, 0.5)
											-- (2.5+2, 0.5);
											
    \draw[dashed, line width=0.4mm] (2.5+1.125+0.04-0.2, 0.5) arc (180:95:0.2);
%    \draw[dashed, line width=0.6mm] (0.625+0.04, 0.5) -- (0.875+0.04, 0.5);
	\draw (2.5+1.02, 0.8) node {$\alpha$};
%	\draw[|-|,thick, line width=0.3mm] (2.5+1.375+0.04, 0.35) -- (2.5+1.625+0.04, 0.35) node[sloped,midway,below] {$W_F$};

% Third geometry
   	\draw[line width=0.6mm] (5+0,0) rectangle (5+2,2);
   	
   	\draw[line width=0.6mm] (5+0, 0.5)	-- (5+0.125+0.08, 0.5)
										-- (5+0.125-0.08, 1.5)
										-- (5+0.375-0.08, 1.5)
										-- (5+0.375+0.08, 0.5)
										-- (5+0.625+0.08, 0.5)
										-- (5+0.625-0.08, 1.5)
										-- (5+0.875-0.08, 1.5)
										-- (5+0.875+0.08, 0.5)
										-- (5+1.125+0.08, 0.5)
										-- (5+1.125-0.08, 1.5)
										-- (5+1.375-0.08, 1.5)
										-- (5+1.375+0.08, 0.5)
										-- (5+1.625+0.08, 0.5)
										-- (5+1.625-0.08, 1.5)
										-- (5+1.875-0.08, 1.5)
										-- (5+1.875+0.08, 0.5)
										-- (5+2, 0.5);
										
    \draw[dashed, line width=0.4mm] (5+1.125+0.08-0.2, 0.5) arc (180:100:0.2);
%    \draw[dashed, line width=0.6mm] (0.625+0.08, 0.5) -- (0.875+0.08, 0.5);	
	\draw (5+1.05, 0.8) node {$\alpha$};
%	\draw[|-|,thick, line width=0.3mm] (5+1.375+0.08, 0.35) -- (5+1.625+0.08, 0.35) node[sloped,midway,below] {$W_F$};

\end{tikzpicture}}
	\vspace{-0.3cm}	
	\caption{Internal morphology with nanorods with a varying incidence angle $\alpha$.}
	\label{fig:angle_geo}
\end{figure}

The considered device geometry is shown in Fig.~\ref{fig:angle_geo}, 
where the number of rods is kept constant to four for each material but the incidence angle $\alpha$ is varied in a range from $90^{\circ}$ to $77^{\circ}\,11'$. The geometric data are $L_{cell} = L_{elec} = 150\,\nano\meter$, $L_R = 75\,\nano\meter$ and $W_R = 37.5\,\nano\meter$.
\begin{figure}[bt!]
	\centering
	\subfigure[Interface length]{
		\includegraphics[width=0.475\textwidth]{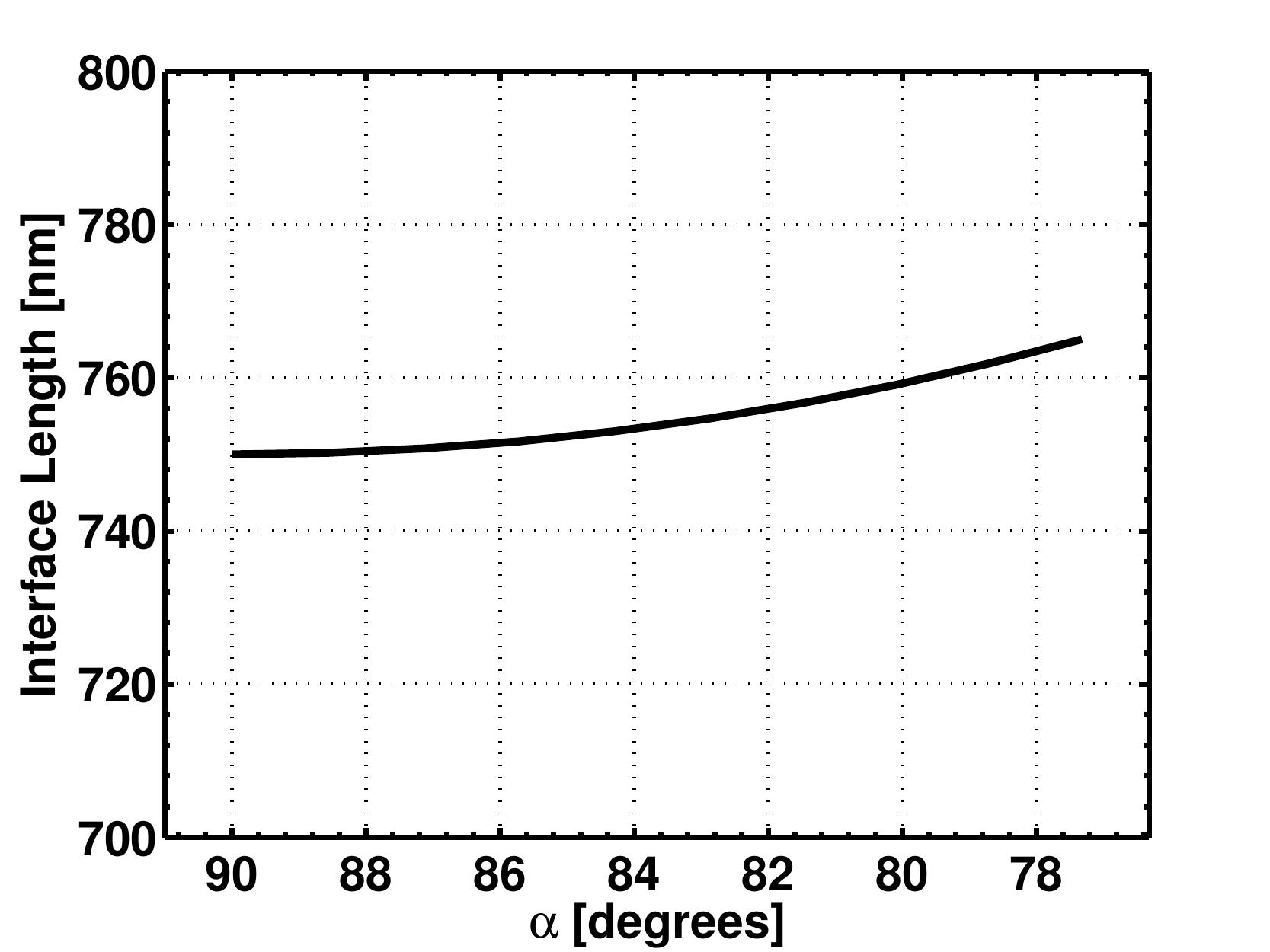}
		\label{fig:angles_left}
	}
	\subfigure[Short circuit current density]{
		\includegraphics[width=0.475\textwidth]{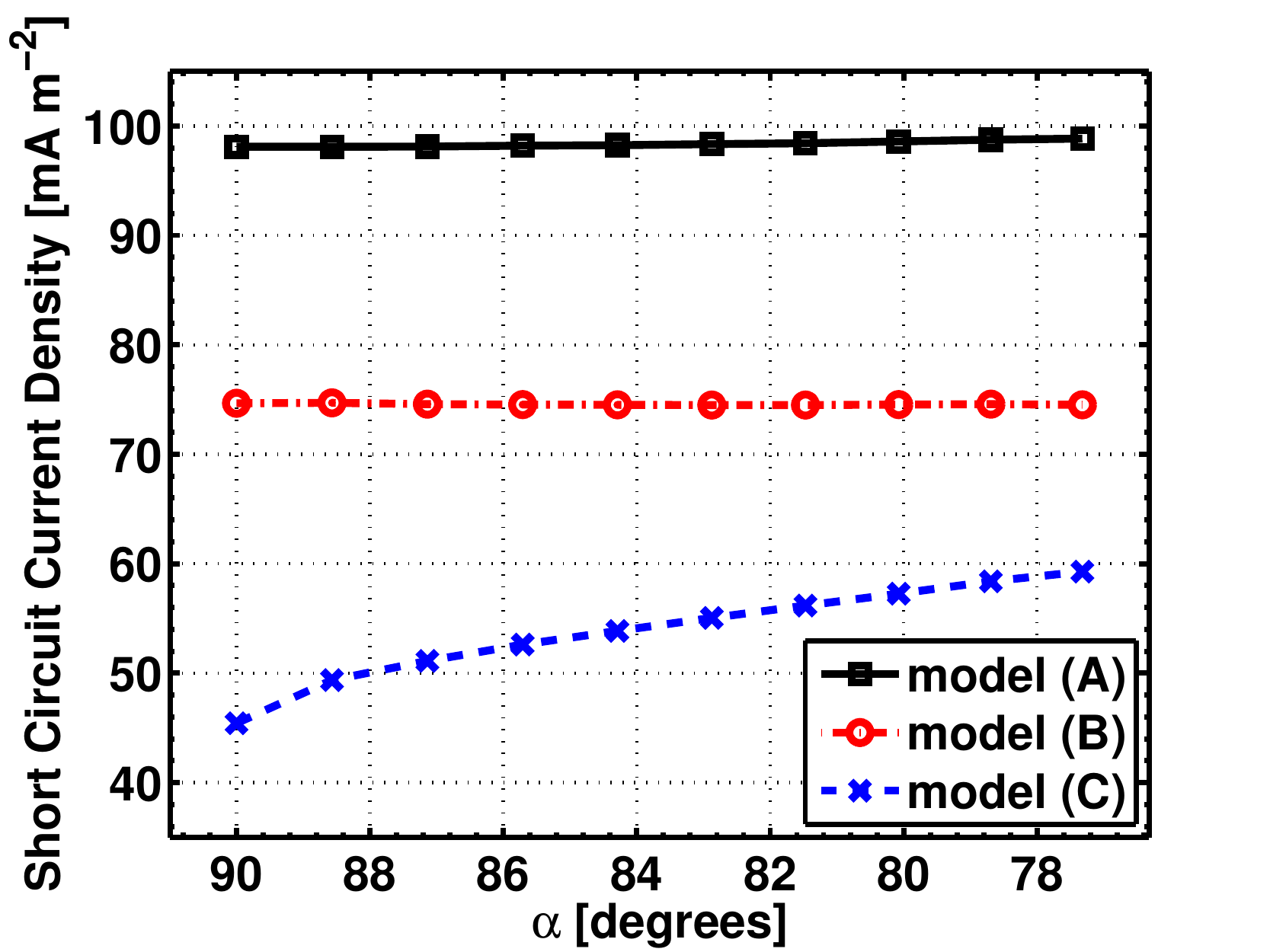}
		\label{fig:angles_right}
	}
	\caption{Interface length and short circuit current density as functions of $\alpha$.}
	\label{fig:angles}
\end{figure}

Since the changes in the amplitude for $\alpha$ are small, the interface length does not vary significantly (as demonstrated by Fig.~\ref{fig:angles_left}) and we expect model~(A) to be quite insensitive to such small modifications since $E_y$ mainly depends on the potential drop
across the electrodes. Concerning with model~(B), we again do not
expect a relevant sensitivity to such variations in the interface morphology 
since the changes in $E_n$ and $E_t$ should balance in the overall contribution. We instead expect
model~(C) to be most sensitive since the normal field that is screened 
at the interface may experience significant variations 
as a function of the angle $\alpha$.

Our expectations are confirmed by the results in Fig.~\ref{fig:angles_right}, showing that the performance of the device in terms of computed short circuit current density does not vary with $\alpha$ when models 
(A) and (B) are considered, while if model (C) 
is used, an increase of the short circuit 
current density is observed as soon as the inclination of the 
nanorod structure is modified with respect to the initial 
configuration. This behavior can be explained as follows. 
The choice of model~\eqref{eq:kdissmodel3} predicts an increase of the dissociation for negative values of the normal electric field that is higher than the reduction for positive values of $E_n$. 
Since at short circuit the electric field can be reasonably assumed to be
directed along the $y$ axis (\emph{i.e.}, from the cathode to the anode), 
the sides of each rod experience opposite normal fields. As a result, the overall effect is dominated by the contribution of the sides with negative fields and dissociation is enhanced.
%
%
%This is due to the facts that in these configurations $E_n$ has opposite sign on the sides of the rods and in the $k_{diss}$ model the increase of the dissociation rate for negative field values is more than the dumping for positive fields. Anyway, this behavior seems to be physically unmotivated and again the version of the model that includes also the tangential component is the one that gives the most reasonable results.

\subsection{The Case of a Complex Interface Morphology}
\label{sec:morphology_2}

\begin{figure}[hbt!]
	\centering
%	\subfigure{\includegraphics[width=0.4\textwidth, angle=90]{figures/our_mesh_small.pdf}}
	\includegraphics[width=0.35\textwidth, angle=90]{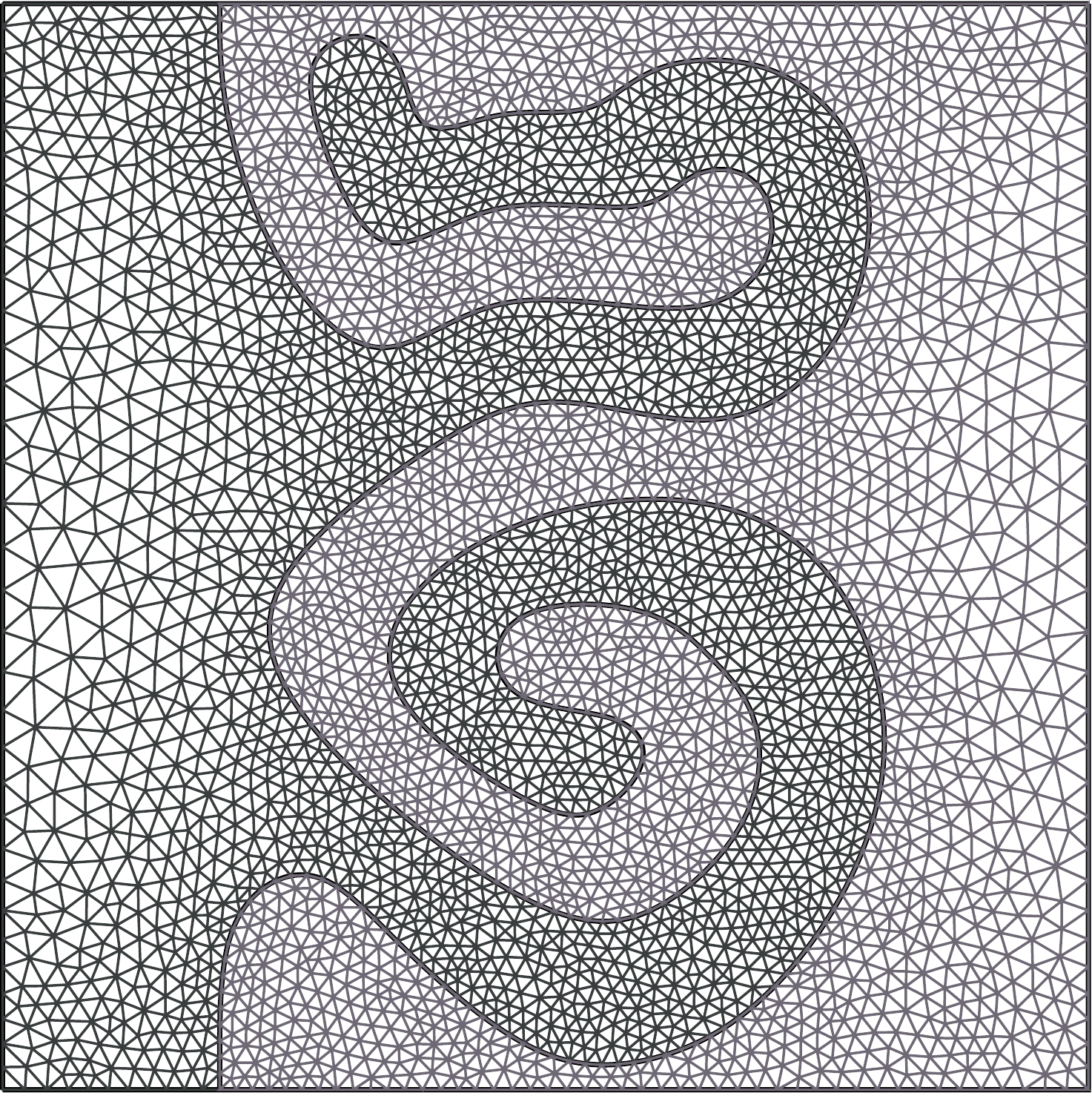}
%	\hspace{0.5cm}
%	\subfigure{\includegraphics[width=0.4\textwidth, angle=90]{figures/their_mesh_small.pdf}}
%	\caption{Comparison between the mesh used to numerically solve
%	the model of Sect.~\ref{sec:reduced_model} and the one that should 
%	be used to numerically solve
%	the full-scale model of Sect.~\ref{sec:fullscale_model}.}
% 	in the simulation and another one to be used in case of a model with volumetric terms in the neighborhood of the interface.}
	\caption{The computational mesh used to numerically solve
	the model of Sect.~\ref{sec:macroscale_model} in the 
	case of a complex geometry.}
	\label{fig:mesh}
\end{figure}

In this concluding section, we test the versatility of the 
model proposed in the present article in dealing with a 
very complex internal morphology as that shown in Fig.~\ref{fig:mesh}.
% aim at demonstrating the benefit
% produced by the model proposed in the present article, compared
% to existing approaches, in dealing with a very complex internal morphology
% as that shown in Fig.~\ref{fig:mesh}.
In this regard, it is important to notice that the use of the microscale 
model~\eqref{eqs:excitons}-\eqref{eq:polarons} would require 
an extremely fine grid resolution to accurately 
describe the volumetric terms in the active layer 
around the donor-acceptor interface, while the 
use of the macroscale model~\eqref{eq:lumped_model} 
has the twofold advantage of
considerably simplifying the design of the computational mesh 
and reducing the size of the nonlinear algebraic system to be solved. 

\begin{figure}[!hbt]
	\centering
	\subfigure[Charge carrier density]{
		\includegraphics[height=0.38\textwidth]{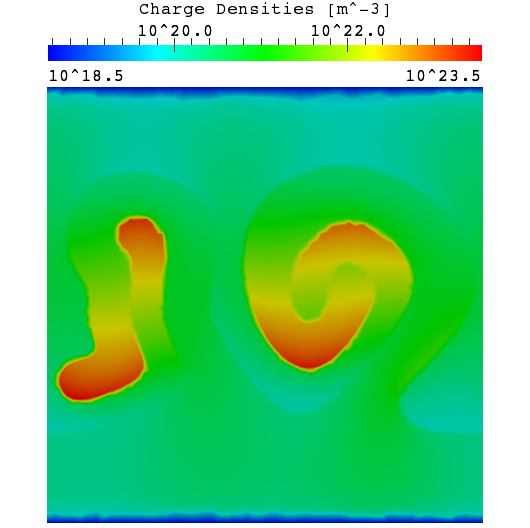}
		\label{fig:complex_geo_left}
	}
	\subfigure[Current-voltage characteristics]{
		\includegraphics[height=0.38\textwidth]{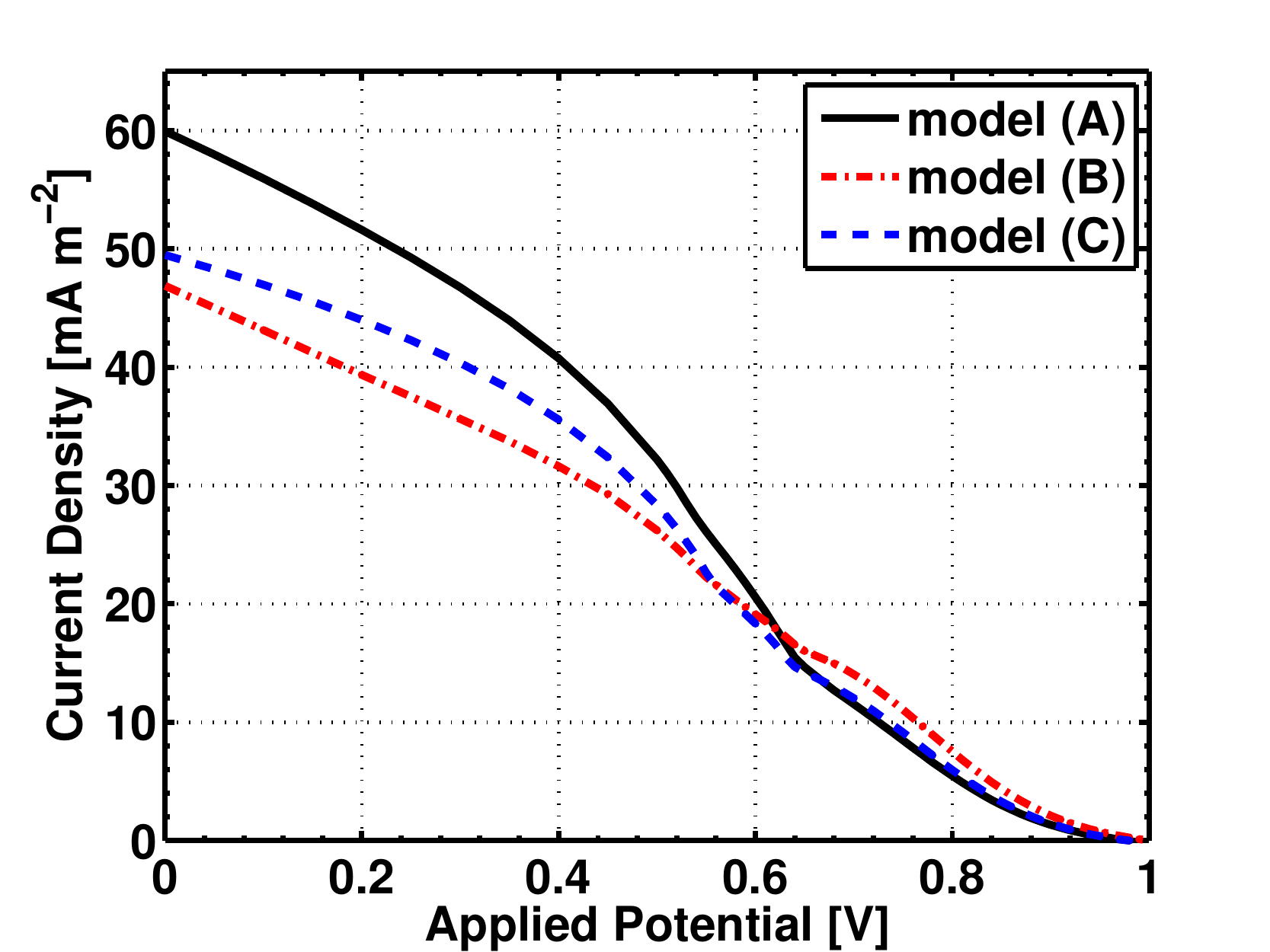}
		\label{fig:complex_geo_right}
	}
	\caption{Log-plot of charge carrier density [$\meter^{-3}$]
	and current-voltage characteristics 
	for a device with very complex internal structure.}
	\label{fig:complex_geo}
\end{figure}

Fig.~\ref{fig:complex_geo_left} illustrates the computed charge carrier density at short circuit ($V_{appl} = 0\,\volt$) for the geometry of Fig.~\ref{fig:mesh}, where the domain is a square $150\,\nano\meter$ sided, with exciton generation rate
$Q = 1.53 \cdot 10^{25}\,\meter^{-3}\,\second^{-1}$ and using 
model (B) for $k_{diss}$.
Notice in particular that the densities assume much higher magnitudes 
compared to those of Fig.~\ref{fig:densities}. This is a consequence of the complexity of the geometry, where donor and acceptor form dead-end areas in which the charges are trapped and experience recombination. In Fig.~\ref{fig:complex_geo_right} we show again a comparison of the current-voltage characteristics obtained using the three different polaron dissociation rate models. The differences among the obtained characteristic lines are reduced with respect to the previous simulated cases. In particular the computed short circuit current densities attain closer values with respect to more regular morphologies,
such as that of Fig.~\ref{fig:finger_geo}, for comparable values of the interface length (approximately $900\,\nano\meter$), see Fig.~\ref{fig:isc_vs_length}.
This is probably to be ascribed to the tortuosity of device internal
morphology which makes interface recombination effects more significant
than in the case of a more regular internal structure.
% 
% 
% DA QUI IN POI FORSE NELLE CONCLUSIONI\\
% Since the effectiveness of this procedure is very good, it could be in principle possible to simulate the behavior of bulk-heterojunction solar cells by directly considering the actual interior morphology without employing an averaged homogeneous media approach as it has been done in \cite{hwang_greenham,defalco}. In this way more details of the cell operation would be available and it would be possible to achieve a better understanding of current photogeneration.

\section{Concluding Remarks and Future Perspectives}
\label{sec:concl}
The research activity object of the present article 
is a continuation of the mathematical study
of organic photovoltaic devices started in~\cite{defalco}
and is focused to:
\begin{description}
\item[-] the accurate and computationally efficient 
modeling of  photoconversion mechanisms occurring at the
interface separating the acceptor and donor layers;
\item[-] 
the investigation of the impact of the interface
morphology and of polaron pair dissociation 
on device performance.
\end{description}

With this aim, we
propose a two-scale (micro- and macro-scale) 
multi-dimensional model for organic solar cell devices with 
arbitrary interface geometries. 
The microscale model is a system of 
incompletely parabolic nonlinear PDEs in drift-diffusion form
set in a heterogeneous domain.
The macroscale model is obtained via a micro-to-macro scale 
transition that consists of averaging the mass balance 
equations in the direction normal to the interface,
giving rise to nonlinear transmission conditions 
parametrized by the interfacial width.
This averaging procedure is similar to model-reduction 
techniques used in porous media 
with thin fractures~\cite{Jaffre2005}, in reaction problems 
with moving reaction fronts~\cite{transitions} and
in electrochemical transport across biological 
membranes~\cite{Mori2009}.
The fact that in the macroscale model the interface is reduced 
to a zero width surface is further exploited to account for the
local dependence of the polaron dissociation rate on 
the electric field orientation, which is the main advantage 
--together with the computational cost reduction-- of our
approach, as compared to previous multi-dimensional 
models~\cite{buxton2007,williams2008,jerome2010}.

%% Our formulation, compared to ,
%% has the advantage of simplifying
%% the discrete computational domain
%% while allowing, at the same time, to easily include 
%% in the model a local dependency of the coefficients 
%% on the orientation of the electric field.
%% This latter feature of the proposed model 
%% has an important consequence 
%% in the optimal design of novel
%% advanced photoconversion devices, as the dissociation
%% rate is a critical physical parameter affecting the
%% efficiency performance of the solar cell.

\null

Extensive numerical simulations of realistic
device structures are carried out to study the performance of the proposed
models and the impact of the lumping procedure. 
First, one-dimensional
transient simulations under different working
conditions are carried out to verify the accuracy 
of the macroscale model with respect to the microscale system. 
Results indicate that in the physically reasonable range of values for the 
parameter $H$ the relative discrepancy between the micro and macroscale formulations
is consistently below 10\%. 
Two-dimensional realistic device structures with various interface
morphologies are then numerically investigated to assess
the impact of our novel model for $k_{diss}$ on the main device properties
(short circuit current and open circuit voltage).
Simulation results indicate that, if the electric field orientation
relative to the interface is taken into due account, the device
performance is determined not only by the total
interface length but also by its shape.

%Two-dimensional realistic device structures are then considered,
%the focus being the investigation of the principal device design parameters 
%(short circuit current and open circuit voltage) 
%as functions of material interface morphology. In particular, a 
%novel approach to consistently account for interface electric field
%orientation in the polaron dissociation rate is proposed and
%successfully compared to previously existing (simplified) 
%expressions depending on a
%suitable average of the electric field.

\null

Research topics currently under scrutiny include:
\begin{description}
\item[-] application of the proposed computational model to the
study of more complex three-dimensional morphologies,
as considered in~\cite{kimber2010};
\item[-] investigation of more advanced 
models for carrier mobilities and polaron dissociation rate,
as well as the simulation of other material blends 
currently employed in the fabrication 
of up-to-date organic solar cells (see, 
{\em e.g.},~\cite{beverina,beverina2,caironi2007});
%\item[-] inclusion of time dependency in the computational
%experiments, in order to extend to bilayer OSCs the 
%analysis carried out in~\cite{defalco};
\item[-] extension of the model to the general case where
$\Gamma^+$ and $\Gamma^-$ are free boundaries to be determined
for each $\vect{x} \in \Gamma$ and at
each time level $t>0$;
% \item[-] \rs{development of suitable adaptive mesh refinement 
% strategies to cope with possible corner singularities arising in
% the solution in the
% case of irregularly shaped interface morphologies;}
\item[-] a more thorough mathematical investigation of 
the proposed equation 
system~\eqref{eq:lumped_model} in both stationary
and time-dependent regimes to extend the
analysis carried out in~\cite{defalco}.

\end{description}

\section{Acknowledgements}
\label{sec:ack}
The authors 
%acknowledge the contribution of 
%Maurizio M.\ Cogliati to the research topic addressed in this article, and 
wish to thank Dr.\ Dario Natali from  
Dipartimento di Elettronica e Informazione, Politecnico di Milano 
and Dr.\ Mos\`e Casalegno from Dipartimento di Chimica 
``Giulio Natta'', Politecnico di Milano for many fruitful 
and stimulating discussions.
The third author was supported by the M.U.R.S.T.\ grant nr.\ 
200834WK7H005 {\em Adattivit\`a Numerica e di Modello
per Problemi alle Derivate Parziali} (2010-2012).

%\section{Appendix}
%\label{sec:appendix}
\appendix
\section{Finite Element Discretization}
\label{sec:appendix}

In this appendix, for sake of completeness, 
we provide more detail about the Finite Element (FE) 
discretization of the linear problem resulting from the application
of time semi-discretization and linearization to the 
equation system~\eqref{eq:lumped_model} as schematically described 
in Sect.\ref{sec:numer_approx}. 
Before we proceed we need to introduce some notation. 

The time semi-discretization consists of approximating the
time derivative of the generic quantity $U$ 
(representing any of the unknowns in~\eqref{eq:lumped_model}) as
\begin{equation}
\dfrac{\partial U}{\partial t} \simeq w_{0} U_{N} + 
\displaystyle\sum_{m = 1}^{m = M} w_{m} U_{N - m} =
w_{0} U_{N} + 
d_{N,M} \left( U \right),
\end{equation}
where $N$ is the index of the current time step and $M$ is 
the order of the adopted BDF formula.
The notation $d_{N,M} \left( U \right)$ is used 
to group together the terms that depend only on 
results from past time steps and is therefore a known 
quantity at the $N$-th time integration level.

To treat the spatial discretization of the problem, we assume only
for ease of presentation that $\Omega$ is a rectangular domain, 
as depicted in Fig.~\ref{fig:Omega_2d}, but the approach remains completely
valid also in the three-dimensional case, provided to replace 
``triangle" by ``tethrahedron" and ``edge" by ``face".
Let $\mathcal{T}_{h}$ denote a conformal partition into open triangles
$K$ of the computational domain $\Omega$, 
$h$ being the maximum diameter over all triangles,  
and let $\Omega_{n,h}$ and $\Omega_{p,h}$ denote the finite 
element partitions of the subregions $\Omega_n$ and $\Omega_p$, such
that $\Gamma_{h} = \partial \Omega_{n,h} \cap \partial 
\Omega_{p,h}$ is their separating interface, 
consisting of the union of a set of edges of $\mathcal{T}_{h}$.

We introduce the following finite dimensional spaces of FE functions:
\begin{subequations}
\label{eq:fem_spaces}
\begin{align}
&V_{h} \equiv \left\{v \in C^{0}\left(\overline{\Omega}\right), \left.v\right\vert_{K} 
\in \mathbb{P}_{1}(K)\, \forall K \in \mathcal{T}_{h}\right\}&\\
&V_{h}^{g} \equiv \left\{v \in V_{h}, v = g \mbox{ at the nodes on }
\Gamma_{A}\cup\Gamma_{C}, g \in  C^0\left( \Gamma_{A}\cup\Gamma_{C}\right)
\right\} & \\
&V_{n, h} \equiv \left\{\left.v\right|_{\Omega_{n,h}}, v \in V_{h} \right\}&\\
&V_{p, h} \equiv \left\{\left.v\right|_{\Omega_{p,h}}, v \in V_{h} \right\}&\\
&V_{\Gamma, h} \equiv \left\{\left.v\right|_{\Gamma_{h}}, v \in V_{h} \right\}.
\end{align}
\end{subequations}
Let $\varphi_D \in C^0\left(\Gamma_{A}\cup\Gamma_{C}\right)$ be such that 
$\left.\varphi_D\right|_{\Gamma_{C}} = 0$ and 
$\left.\varphi_D\right|_{\Gamma_{A}} = V_{appl} + V_{bi}$.
Then, we denote by
\begin{equation}
\mathbf{y}_{h} = [e_{h}, \widetilde{P}_{h}, n_{h}, p_{h}, \varphi_{h}]^{T} 
\in 
\left(V_{h}^{0} \times V_{\Gamma, h} \times V_{n, h} \times V_{p, h} \times V_{h}^{\varphi_D}\right) \equiv  \mathbf{Y}_h
\end{equation}
the vector of discrete unknown functions at a given quasi-Newton iteration of a given time step, and by 
\begin{equation}
\label{eq:fem_increments}
\delta \mathbf{y}_{h} = [\delta e_{h}, \delta \widetilde{P}_{h}, 
\delta n_{h}, \delta p_{h}, \delta \varphi_{h}]^{T}
\in \left(V_{h}^{0} \times V_{\Gamma, h} \times V_{n, h} 
\times V_{p, h} \times V_{h}^{0}\right) \equiv \mathbf{V}_h
\end{equation}
the corresponding increments to be computed in order to advance to the next iteration of the quasi-Newton method.

The linear problem to be solved in order to compute the increments~\eqref{eq:fem_increments} reads: \\
given $\mathbf{y}_{h} \in \mathbf{Y}_h$, 
find $\delta \mathbf{y}_{h} \in \mathbf{V}_h$ such that:
\begin{subequations}
\label{eq:fem_problem}

\begin{align}
&
\displaystyle\int_{\Omega} D_{e} \nabla \delta e_{h} \cdot \nabla v +
\int_{\Omega} \left( \dfrac{1}{\tau_{e}} + w_{0}\right) \delta e_{h} \ v +
\int_{\Gamma_{h}} \left[ \dfrac{2H}{\tau _{diss}}\ \delta e_{h} -
 \eta k_{rec}\ \delta \widetilde{P}_{h} \right]\ v = &\notag \\ \label{e_red_fem} &
\qquad -\left[\displaystyle\int_{\Omega} D_{e} \nabla e_{h} \cdot \nabla v
+ \int_{\Omega} f_{e} (\mathbf{y}_{h})\ v\ + \int_{\Gamma_{h}}\ g_{e} (\mathbf{y}_{h})\ v\right],
& \qquad \forall v \in V_{h}^{0} &
\\[3mm]
\label{Polaron_red_fem}
&
\int_{\Gamma_{h}} \left[ -\dfrac{2H}{\tau _{diss}}\ \delta e_{h} +  
\left( w_{0} + k_{diss} + k_{rec} \right) \delta \widetilde{P}_{h} -
2H\gamma\ \left( p_{h} \delta n_{h} + n_{h} \delta p_{h} \right) \right]\ v = 
- \int_{\Gamma_{h}}\ g_{\widetilde{P}} (\mathbf{y}_{h})\ v,&
\qquad \forall v \in V_{\Gamma, h}, &
\\[3mm]
&
\displaystyle\int_{\Omega_{n,h}} 
\left( \mathbf{D}_{n,h} \nabla \delta n_{h} - \boldsymbol{\mu}_{n,h} \delta n_{h} \nabla \varphi_{h}
\right) \cdot \nabla v 
- \displaystyle \int_{\Omega_{n,h}} 
\mu _{n} n_{h} \nabla \delta \varphi_{h} \cdot \nabla v  + 
\displaystyle\int_{\Omega_{n,h}} w_{0}\ \delta n_{h}\ v + & \notag \\ 
& \qquad
\phantom{- }\int_{\Gamma_{h}} \left[ k_{diss} \delta \widetilde{P}_{h} + 
2H \gamma \left( p_{h}\delta n_{h} + n_{h}\delta p_{h}  \right)\right] v + 
\int_{\Gamma_{C}} \dfrac{\alpha_n}{\kappa_n} \delta n_{h}\ v
= & \notag \\ \label{n_red_fem}
& \qquad
- \left[\displaystyle\int_{\Omega_{n,h}} \left( \mathbf{D}_{n,h} \nabla n_{h} - 
\boldsymbol{\mu}_{n,h} n_{h} \nabla \varphi_{h} \right) \cdot \nabla v
+ \int_{\Omega_{n,h}} f_{n} (\mathbf{y}_{h})\ v\ 
+ \int_{\Gamma_{h}}\ g_{n} (\mathbf{y}_{h})\ v +
  \int_{\Gamma_{C}} b_{n}\left(\mathbf{y}_{h}\right) \ v
	\right],
&
\qquad \forall v \in V_{n, h}, &
\\[3mm]
&
\displaystyle\int_{\Omega_{p,h}} 
\left(
\mathbf{D}_{p,h} \nabla \delta p_{h} + \boldsymbol{\mu}_{p,h} \delta p_{h} \nabla \varphi_{h} 
\right) \cdot \nabla v 
+ \displaystyle\int_{\Omega_{p,h}} 
\mu _{p} p_{h} \nabla \delta \varphi_{h} \cdot \nabla v + 
\displaystyle\int_{\Omega_{p,h}} w_{0}\ \delta p_{h}\ v + & \notag \\ 
& \qquad
\int_{\Gamma_{h}} \left[ k_{diss} \delta \widetilde{P}_{h} + 
2H \gamma \left( p_{h}\delta n_{h} + n_{h}\delta p_{h}  \right)\right] v + 
\int_{\Gamma_{A}} \dfrac{\alpha_p}{\kappa_p} p_{h}\ v
=& \notag \\ \label{p_red_fem}
& \qquad
- \left[\displaystyle\int_{\Omega_{p,h}} \left( \mathbf{D}_{p,h} \nabla n_{h} + 
\boldsymbol{\mu}_{p,h} p_{h} \nabla \varphi_{h} \right) \cdot \nabla v
+ \int_{\Omega_{p,h}} f_{p} (\mathbf{y}_{h})\ v\ +
\int_{\Gamma_{h}}\ g_{p} (\mathbf{y}_{h})\ v +
\int_{\Gamma_{A}} b_{p}\left(\mathbf{y}_{h}\right) \ v\right],
 &
 \qquad \forall v \in V_{p, h}, &
\\[3mm]
\label{Poisson_red_fem}
&
\displaystyle\int_{\Omega} 
\varepsilon \nabla \delta \varphi_{h} 
\cdot \nabla v + 
q \int_{\Omega_{n,h}} \delta n_{h} -
q \int_{\Omega_{p,h}} \delta p_{h}
 = - \left[\displaystyle\int_{\Omega} 
\varepsilon \nabla \varphi_{h} 
\cdot \nabla v 
+ \int_{\Omega} f_{\varphi} (\mathbf{y}_{h})\ v \right]&
\qquad \forall v \in V_{h}^{0}, &
\end{align}
\end{subequations}
where  
$v$ denotes in each of the equations~\eqref{eq:fem_problem} 
the test function in the appropriate FE space, 
and where we have made use of the following definitions: 
\begin{subequations}
\begin{align}
&f_{e} \left( \mathbf{y}_{h} \right) = \left( \dfrac{1}{\tau_{e}} + w_{0} \right) e_{h} - Q  + d_{N,M} \left( e_{h} \right)&\\
&g_{e} \left( \mathbf{y}_{h} \right) = \dfrac{2H}{\tau_{diss}} e_{h} - \eta k_{rec} \widetilde {P}_{h}&\\
&g_{\widetilde{P}} \left( \mathbf{y}_{h} \right) =
- \dfrac{2H}{\tau_{diss}} e_{h} + \left(w_{0} + k_{diss} + k_{rec} \right) \widetilde{P}_{h} 
- 2H\gamma n_{h} p_{h} + d_{N,M} \left( \widetilde{P}_{h} \right)&\\
&f_{n} \left( \mathbf{y}_{h} \right) =  w_{0} n_{h} + 
d_{N,M} \left( n_{h} \right) &\\
&g_{n} \left( \mathbf{y}_{h} \right) = 
k_{diss} \widetilde{P}_{h} + 2H\gamma n_{h} p_{h}&\\
&b_{n} \left( \mathbf{y}_{h} \right) = \dfrac{1}{\kappa_{n}} \left\{ \alpha_{n} n_{h} - \beta_{n} \right\}&\\
&f_{p} \left( \mathbf{y}_{h} \right) = w_{0} p_{h} 
+ d_{N,M} \left( p_{h} \right) &\\
&g_{p} \left( \mathbf{y}_{h} \right) = k_{diss} \widetilde{P}_{h} + 2H\gamma n_{h} p_{h}&\\
&b_{p} \left( \mathbf{y}_{h} \right) = \dfrac{1}{\kappa_{p}} \left\{ \alpha_{p} p_{h} - \beta_{p} \right\}&\\
&f_{\varphi} \left( \mathbf{y}_{h} \right) = q n_{h} - q p_{h}.&
\end{align}
\end{subequations}
Then, once system~\eqref{eq:fem_problem} is solved, the unknown vector 
is updated as
$$
\mathbf{y}_{h} \leftarrow  \mathbf{y}_{h} + \delta \mathbf{y}_{h}.
$$
A couple of further comments is in order about the numerically 
stable implementation of the FE linear system~\eqref{eq:fem_problem}.
First, note that $\mathbf{D}_{n,h}$ and $\boldsymbol{\mu}_{n,h}$ 
in~\eqref{n_red_fem} 
(respectively $\mathbf{D}_{p,h}$ and $\boldsymbol{\mu}_{p,h}$ 
in~\eqref{p_red_fem}) are tensor diffusivities and mobilities
chosen according to the Exponential Fitting stabilization technique as in~\cite{sharmacarey89,bank1998,gatti1998,xu1999,lazarov} in order
to avoid the onset of possible spurious oscillations in the 
discrete electron and hole densities due to drift terms.
Second, all the integrals involving zeroth order terms 
are computed using the two-dimensional trapezoidal quadrature rule
in order to end up with strictly positive diagonal (approximate)
mass matrices~\cite{bankrose87}.

\bibliographystyle{amsplain}
\bibliography{paper_osc_multiscale}

\end{document}